\documentclass[final,5p,times,twocolumn,fleqn]{elsarticle}

\usepackage{amsmath}
\usepackage{amssymb,amsthm, algorithm, algpseudocode,xcolor}
\usepackage[utf8]{inputenc}
\usepackage[hidelinks]{hyperref}
\usepackage{multicol}
\usepackage{wrapfig}
\usepackage{float}
\usepackage{multirow}
\usepackage{cleveref} 
\usepackage{booktabs}  
\usepackage{array}     

\journal{Neural Network}

\begin{document}
\newcommand{\customcite}[2]{[\cite{#1}, Section #2]}
\newtheorem{lemma}{Lemma}[section] 
\newtheorem{proposition}{Proposition}[section] 
\newtheorem{definition}{Definition}[section]
\newtheorem{theorem}{Theorem}[section]
\newtheorem{assumption}{Assumption}[section]
\begin{frontmatter}



\title{FxTS-Net: Fixed-Time Stable Learning Framework for Neural ODEs}

\author[1]{Chaoyang Luo}

\author[2]{Yan Zou}

\author[1]{Wanying Li}

\author[1]{Nanjing Huang \corref{cor1}}
\ead{nanjinghuang@hotmail.com}\ead{njhuang@scu.edu.cn}

\cortext[cor1]{Corresponding author}

\affiliation[1]{organization={Department of Mathematics},
	addressline={Sichuan University},
	city={Chengdu},
	postcode={Sichuan, 610064},
	country={China}}
\affiliation[2]{organization={Department of Artificial Intelligence and Computer Science},
	addressline={Yibin University},
	city={Yibin},
	postcode={Sichuan, 644000},
	country={China}}
	
\begin{abstract}
%
Neural Ordinary Differential Equations (Neural ODEs), as a novel category of modeling big data methods, cleverly link traditional neural networks and dynamical systems. However, it is challenging to ensure the dynamics system reaches a correctly predicted state within a user-defined fixed time. To address this problem, we propose a new method for training Neural ODEs using fixed-time stability (FxTS) Lyapunov conditions. Our framework, called FxTS-Net, is based on the novel FxTS loss (FxTS-Loss) designed on Lyapunov functions, which aims to encourage convergence to accurate predictions in a user-defined fixed time. We also provide an innovative approach for constructing Lyapunov functions to meet various tasks and network architecture requirements, achieved by leveraging supervised information during training. By developing a more precise time upper bound estimation for bounded non-vanishingly perturbed systems, we demonstrate that minimizing FxTS-Loss not only guarantees FxTS behavior of the dynamics but also input perturbation robustness. For optimising FxTS-Loss, we also propose a learning algorithm, in which the simulated perturbation sampling method can capture sample points in critical regions to approximate FxTS-Loss. Experimentally, we find that FxTS-Net provides better prediction performance and better robustness under input perturbation.

\end{abstract}


\end{frontmatter}


\section{Introduction}
\label{sec:Introduction}

Neural ODEs have emerged as a promising area of research, offering significant advancements in modeling complex, large-scale datasets through continuous dynamical systems. Originating from the reinterpretation of residual networks (ResNets) as continuous-time dynamical systems \cite{Haber_2018, Ruthotto2020}, Neural ODEs extend to model continuous dynamics and address the limitations of discrete architectures \cite{pmlr-v139-kidger21a}. By incorporating techniques from dynamical systems theory, Neural ODEs enable continuous modeling of data and show broad applicability across various domains, including irregular time series modeling \cite{oh2024stable,NEURIPS2020_4a5876b4}, generative modeling \cite{NEURIPS2019_99a40143,pmlr-v139-kidger21b}, wind speed prediction \cite{YE2022118}, and traffic flow forecasting \cite{CHU2024106549}.

However, Neural ODEs  still face some limitations, particularly in ensuring the stability of the learned  system. The standard learning approach, which relies on differentiating through the ODE solution using techniques like the adjoint method \cite{NEURIPS2018_69386f6b,pmlr-v139-kidger21a}, struggles to guarantee stable behavior (informally, the tendency of the system to remain within some invariant bounded set).  This instability leads to slow convergence and inaccurate predictions \cite{pmlr-v162-rodriguez22a}. Even more concerning, a dynamical system lacking stability guarantees can result in significant output distortions when subjected to small input perturbations.
\begin{figure}[t]
	\centering
	\includegraphics[width=\columnwidth]{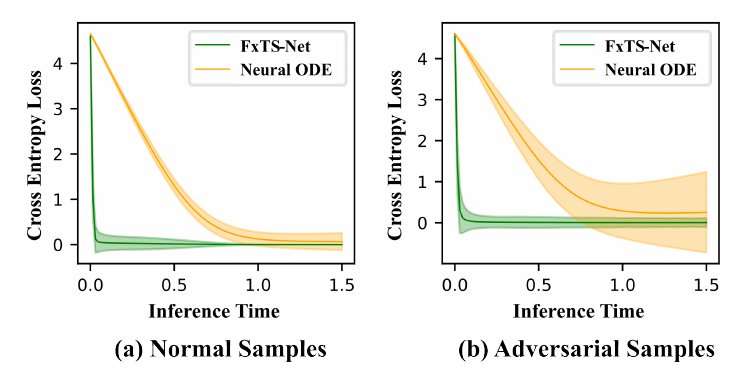}
	\vspace{-20pt}
	\caption{Comparing learned dynamics (plotting inference time vs prediction loss) on 1000 test examples from CIFAR-100. The orange line denotes standard Neural ODEs, the green line denotes FxTS-Net, and the corresponding light-colored areas denote the fluctuating regions across the samples. (a): Experiments on normal samples show FxTS-Net has a faster convergence rate and more stable evolutionary behavior compared to Neural ODEs. (b)  Experiments on adversarial samples indicate that Neural ODEs without fixed-time stabilization result in significant distortions in even small input perturbations.}
	\vspace{-15pt}
	\label{fig:contrast_plot}
\end{figure}

Recently, some research has focused on enhancing the stability of Neural ODEs by applying Lyapunov stability theory \cite{pmlr-v162-rodriguez22a,NEURIPS2019_0a4bbced}. Kolter et al. \cite{NEURIPS2019_0a4bbced} proposed a joint learning approach that integrates Neural ODEs with learnable Lyapunov functions to ensure system stability. However, this approach intensifies the trade-off between stability and accuracy \cite{pmlr-v162-rodriguez22a}. To alleviate the trade-off, Rodriguez et al. introduced a method for training Neural ODEs using Lyapunov exponential stabilization conditions  \cite{pmlr-v162-rodriguez22a}. Since this method relies on the fact that the output function and loss function must satisfy Lyapunov functions conditions, it hampers the applicability in many real-world application scenarios. Constructing appropriate Lyapunov functions remains a longstanding challenge in dynamical systems and control. Therefore, it is essential to investigate and enrich the techniques for constructing Lyapunov functions further to satisfy the requirements of complex network structures and diverse tasks in practical applications.

Despite recent progress in enhancing Neural ODEs stability, most of the existing methods \cite{pmlr-v162-rodriguez22a,NEURIPS2019_0a4bbced,NEURIPS2021_7d5430cf,NEURIPS2023_0a443a00} focus on exponential convergence, neglecting the critical aspect in fixed-time stability. Fixed-time stability is essential in many real-world applications where a system must reliably reach a desired state within a user-defined evolution time. For instance, Cui et al. employed the extensive experimental analysis method to illustrate the critical role of ensuring that a suitable state is reached within a pre-determined time \cite{CUI2023576}; and further Chu et al. extended the evolutionary time to ensure that the dynamical system reaches the desired state under input perturbations, and improved robustness against adversarial samples \cite{https://doi.org/10.1049/cvi2.12248}.  However, to the best of authors' knowledge, Neural ODEs stability within a pre-determined time has not been considered in the literature. The present paper is thus devoted to the fixed-time stability in Neural ODEs under some mild conditons. 

To illustrate the research motivation behind FxTS-Net, we compare FxTS-Net and Neural ODEs in Figure \ref{fig:contrast_plot}. From the orange lines and regions, it indicates that the dynamics of the Neural ODE across the state space do not guarantee stable convergence to correct predictions within a user-defined fixed time. Moreover, input perturbations exacerbate this instability leading to worse convergence speed and prediction accuracy.

In this work, we investigate how to use Lyapunov stability theory to learn fixed-time stable dynamical systems for ensuring the efficientibility and the robustness. Our main contributions can be summerized as follows:

\begin{itemize}
	\item To ensure that the learned Neural ODE achieves correct prediction at a user-defined fixed time, we propose FxTS-Net, which enforces fixed-time stable dynamic inference by introducing FxTS-Loss.
	
	\item  For various tasks and network structures, we propose a method to construct Lyapunov functions, i.e., embedding sub-optimisation problems in training to exploit supervised information for constructing it.
	
	\item For optimising the FxTS-Loss, we propose a learning algorithm, including the simulated perturbation sampling method, i.e., simulating the input perturbation by perturbing extracted features to generate the post-perturbation integral trajectory, which captures the sampling points of the critical region to approximate the FxTS-Loss.
	
	\item Theoretically, we show that minimizing the FxTS-Loss will 1) guarantee that the dynamical system satisfies the FxTS Lyapunov condition; 2) enable robustness of input perturbations by developing a more precise time upper bound estimation for bounded non-decreasing perturbed dynamics.
	
	\item Experimentally, FxTS-Net is competitive or superior in prediction accuracy and improves robustness against various perturbations.  Moreover the experimental results aligns with the theoretical insights, emphasizing the beneficial impact of fixed-time stability on  the overall performance of FxTS-Net .
	
\end{itemize}

\section{Background and Related Work}
\subsection{Supervised Learning of Neural ODEs}

\noindent  \textbf{Neural ODEs.} Following previous research work \cite{pmlr-v162-rodriguez22a}, we will consider a class of data-controlled Neural ODEs. Let $(x,y)$ denote the input data pair where $x$ belongs to $\mathbb{R}^{d_x}$, the mathematical model is
\begin{align}
	 &x_c=\phi \left( {x;{\theta _\phi }} \right),  \label{eq:2.1} \\
		&h(t) = \int_0^t f (s,x_c,h(s);{\theta _f}){\mkern 1mu} ds, \label{eq:2.2}\\
	&\hat{y}_{x}\left(t \right) = \psi \left( {h\left( t \right);{\theta _\psi }} \right),  \label{eq:2.3}
\end{align}
where $ \psi : \mathbb{R}^{d_h} \to \mathbb{R}^{d_{\psi}} $ is an output function with parameter $ \theta_{\psi}$, likewise $\phi$ is an input function, and let  $\theta=(\theta_\phi, \theta_f, \theta_\psi) \in \Theta \subset \mathbb{R}^{d_{\Theta}}$. The idea of the Neural ODE is to use a neural network to parameterize a differential equation that governs the hidden states $h(t) \in H \subset \mathbb{R}^{d_h}$ with respect to time t.

In practice, we set that Equation \eqref{eq:2.2} evolves over the time interval $[0, 1]$, i.e., the upper limit of the integral $t=1$, although theoretically it is possible to choose any time to predict $t$.  Additionally, we assume that the state space $H$ is both bounded and path-connected.

\noindent  \textbf{Model Assumption.} We impose a Lipschitz continuity assumption on neural networks used in Neural ODEs as follows.

\begin{assumption} \label{as:1}
	For any neural network \( F(t, x; \theta_F) : \mathbb{R}_+ \times \mathbb{R}^{d_x} \to \mathbb{R}^{d_F} \) with parameter \( \theta_F \), there exists a positive constant \( L_F > 0 \) such that, for all \( t \geq 0 \) and \( x, x' \in \mathbb{R}^{d_x} \),
	\begin{equation}
		|F(t, x; \theta_F) - F(t, x'; \theta_F)| \leq L_F|x - x'|.
	\end{equation}
\end{assumption}

This assumption is not overly onerous since neural networks with activation functions such as tanh, ReLU, and sigmoid functions generally satisfy the Lipschitz continuity condition in Assumption \ref{as:1} \cite{oh2024stable,NEURIPS2019_95e1533e,Latorre2020Lipschitz}.

\begin{figure*}[h] 
	\centering
	\includegraphics[width=\textwidth]{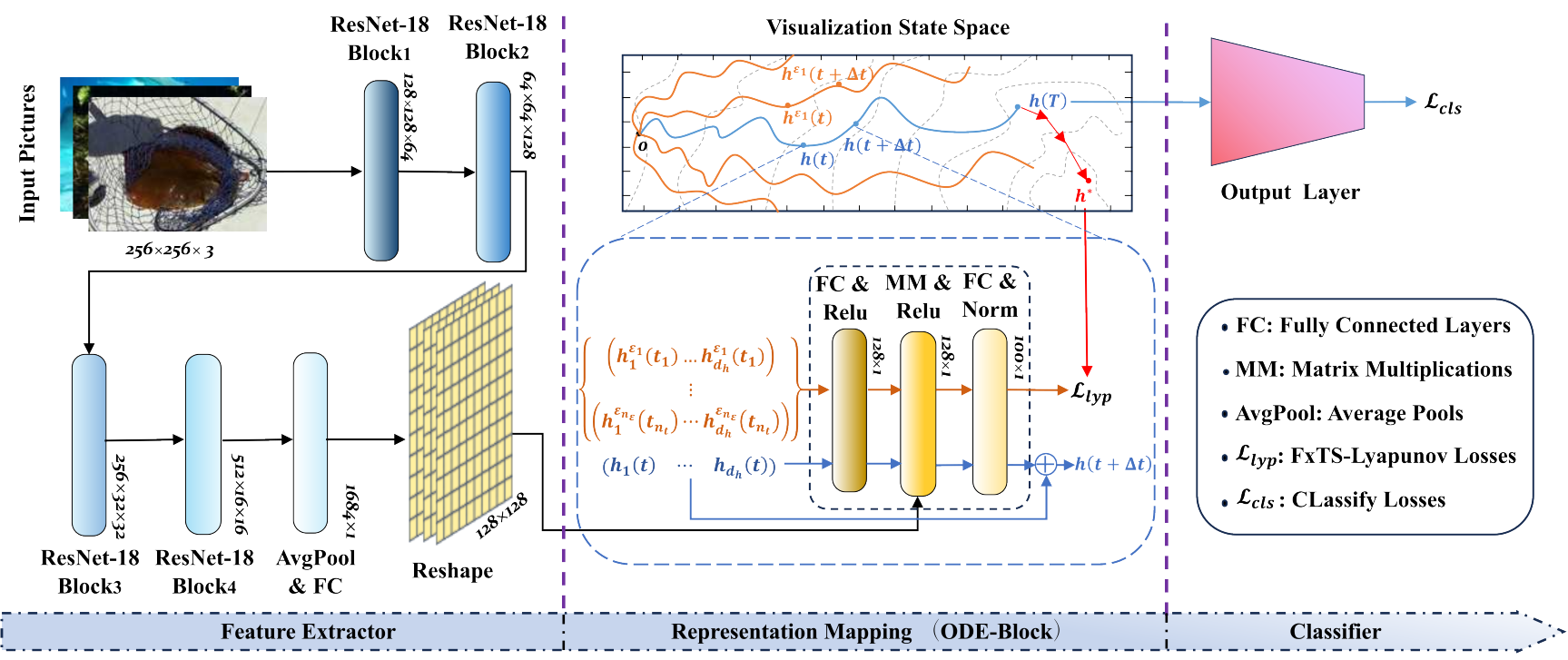} 
	\caption{The overall architecture of FxTS-Net: Using ResNet18 as the feature extractor, a data-controlled Neural ODE as the ODE-Block, a classifier as the output layer, and the corresponding model structure and computational data flow are plotted. Furthermore, on top of the ODE-Block, we plot a phase space, along with the level set of the Lyapunov function $V$, which in this example is minimized  at $h^{*}$. In the phase space, the blue curve indicates the standard integration trajectory, the orange indicates the integration trajectory after perturbing the extracted features, and the red arrow indicates the descent direction of the optimization  problem in Equation \eqref{eq:h}.} 
	\label{fig:fullwidth-top}
\end{figure*}

\noindent  \textbf{Supervised Learning.} Consider the standard supervised learning setup, where the training set of input-output pairs is denoted as  as $(x, y) \sim D$, the optimization model is 
\begin{equation}
\mathop{\arg \min}_{\theta \in \Theta} \sum_{(x,y) \sim D} 	\mathcal{L}(\hat{y}_x(t), y) \label{eq:2.5},
\end{equation}
where $\hat{y}_x{(t)}$ is shorthand for Equations \eqref{eq:2.1}-\eqref{eq:2.3}. The learning goal is to find a parameterization of our model by minimizing a supervised loss over the training data.

Like typical methods in deep learning, the standard approach to training Neural ODEs is backpropagation via Equation \eqref{eq:2.5}. Moreover, the end-to-end training optimization problem can be viewed as the following optimal control problem (for brevity, only using a single $(x,y)$), is given by
\begin{align}
	\mathop{\arg \min}_{\theta \in \Theta} &\sum_{(x,y) \sim D} 	\mathcal{L}(\hat{y}_x(t), y),  \label{eq:Op1} \\
\text{s.t. }  &x_c=\phi \left( {x;{\theta _\phi }} \right), \nonumber \\
&h(t) = h(0) + \int_0^t f (s,x_c,h(s);{\theta _f}){\mkern 1mu} ds, \nonumber\\
&\hat{y}_{x}\left(t \right) = \psi \left( {h\left( t \right);{\theta _\psi }} \right). \nonumber
\end{align}

It is possible to optimize Equation \eqref{eq:Op1}, which can be rolled out of the dynamics using solver backpropagation or concomitant methods \cite{NEURIPS2018_69386f6b,pmlr-v139-kidger21a}. However, in Equation \eqref{eq:Op1}, there is no explicit penalty or regularization to ensure the dynamical system reaches a steady state in a user-defined fixed time. Indeed, we can observe such a problem from Figure \ref{fig:contrast_plot},  where the dynamics of Neural ODEs are not guaranteed to stabilize in a fixed time, leading the Neural ODE to learn a fragile solution. In this work, we propose the FxTS-Net approach to address these limitations through a fixed-time control theory learning objective.

\subsection{Lyapunov Conditions for Fixed-Time Stability} \label{sec:2.2}
Intuitively, a fixed-time stable dynamical system means that all solutions in a region around an equilibrium point flow to that point at a prescribed fixed time. The area of Lyapunov theory \cite{Khalil:1173048} establishes the connection between the various stability and descent according to a particular type of function known as a Lyapunov function. In this work, we assume that the equilibrium point is $h^{*}$.
\begin{definition}[Lyapunov Function]
	A continuously differentiable function \( V: H \rightarrow \mathbb{R} \) is a Lyapunov function if $V$ is positive definite function, i.e., $V(h) > 0$ for $h \neq h^*$ and $V( h^*) = 0$.
\end{definition}
We now recall the fixed-time stability with respect to \(V\). The fixed-time stability (FxTS), defined as the following, allows the settling time to remain uniformly bounded for all initial conditions.
\begin{definition}[FxTS]
	We say that the ODE defined in Equations \eqref{eq:2.1}-\eqref{eq:2.3} is fixed-time stable if there exists a Lyapunov function $V$ and a constant $T_{max} >0 $, such that all solution trajectories $h\left( t\right) $ of the ODE satisfy 
	\begin{equation} \label{eq:exponential_stability}
		  \sup \{ T \geq 0 : V(h(t)) = 0 \text{ for all } t \geq T \} \le T_{max}.
	\end{equation}
\end{definition}

For the study of fixed-time stability, it is necessary to introduce the sufficient
condition for FxTS of the equilibrium point $h^{*}$, which are used to guarantee the convergence of the system trajectory to the equilibrium point in a user-defined fixed time  \cite{GARG2022110314}, as stated in the following theorem.

\begin{lemma}[FxTS Condition\cite{GARG2022110314}] \label{lem:2.1}
	For the ODE in Equations \eqref{eq:2.1}-\eqref{eq:2.3} and a continuously differentiable Lyapunov function \( V \), there exist parameters \(\alpha_1, \alpha_2 > 0\), \(\gamma_1 = 1 + \frac{1}{\mu}\), and \(\gamma_2 = 1 - \frac{1}{\mu}\) with \(\mu > 1\) such that
	\begin{equation} \label{eq:fxs-clf}
		\min_{\theta \in \Theta} \left[ \frac{\partial V}{\partial h}f(h, x_c, t;\theta_{f}) + \alpha_1 V(h)^{\gamma_1} + \alpha_2 V(h)^{\gamma_2} \right]  \leq  0
	\end{equation}
	holds for all \( h \in H \) and \( t \in [0, 1] \). Then there is \( \theta \in \Theta \) that can achieve
	\begin{equation} \label{eq:FxS-stability}
	\frac{\partial V}{\partial h}f(h, x_c, t;\theta_{f})\leq - \alpha_1 V(h)^{\gamma_1} - \alpha_2 V(hx_c)^{\gamma_2},
	\end{equation}
	 and the time of convergence \( T \) satisfies \( T \leq  \frac{\mu \pi}{2\sqrt{\alpha_1 \alpha_2}} \).
\end{lemma}

The introduction of FxTS in Neural ODEs is desirable because: 1) it guarantees fast convergence to desired states (as defined by \( V \)) after integrating for a fixed time (e.g., for \( t \in [0, 1] \)); 2) it has implications for adversarial robustness. Specifically, we require the learned Neural ODE to satisfy the additional structure specified by Equations \eqref{eq:fxs-clf} and \eqref{eq:FxS-stability} for ensuring the fixed-time stability of the Neural ODE. In Section 3, we will develop a learning framework to find a desired parameter \( \theta \) satisfying Equation \eqref{eq:FxS-stability}.

\subsection{Learning Stable Dynamics}

In recent years, many researchers have conducted in-depth studies on how to impose stability on Neural ODEs with various formulations. Using regularization of the flow on perturbed data based on time-invariance and steady-state conditions, Yan et al. proposed TisODE, which outperforms ordinary Neural ODEs in terms of robustness \cite{YAN2020On}. By introducing contraction theory, Zakwan et al.  improved the robustness of Neural ODEs to Gaussian and impulse noise on the MNIST dataset \cite{9809979}. In addition, Kang et al. imposed Lyapunov stability guarantees around the equilibrium point to eliminate the effects of perturbations in the input \cite{NEURIPS2021_7d5430cf}.

We observe that the above works are devoted to the study of the exponential stability for learning dynamical systems. Generally, setting the final evolution time as a pre-defined fixed one is a common requirement in practice and so 
it is crucial to ensure the stationarity of the dynamical system solution within a user-defined fixed time. In this work, we propose a fixed-time stabilization learning framework that allows the learned dynamical system to be stabilized within a pre-defined fixed time.

\section{Fixed-Time Stable Learning Framework} \label{sec:3.1}

The intuition behind the proposed approach is straightforward: as mentioned in Section \ref{sec:2.2}, our goal is to find the parameters \( (\theta_\phi, \theta_f, \theta_\psi) \) of the Neural ODEs that satisfy the FxTS Lyapunov stability condition (Lemma \ref{lem:2.1}). We formulate this approach in two steps: 
\begin{enumerate}
	\item In Section 3.1, for any given output and loss functions, we construct an appropriate Lyapunov function \( V \) based on the supervised information.
	\item In Section 3.2, we design a new FxTS-Loss via the Lyapunov function \( V \), which quantifies the degree of violation for the FxTS contraction condition in Equation \eqref{eq:FxS-stability}. By optimizing FxTS-Loss, the Neural ODE has a stable evolutionary behavior within a pre-defined fixed time.
\end{enumerate}

Theoretically, by optimizing FxTS-Loss, we show that the learned Neural ODE stabilizes to the optimal state point such that the supervised loss is minimized at a pre-defined fixed time (Theorem \ref{th:1}), and obtain a new adversarial robustness guarantee (Theorem \ref{the:2}). The approach architecture can be seen in Figure \ref{fig:fullwidth-top}, and the corresponding learning algorithm is described in Section \ref{sec:learning_algorithms}.

\subsection{ Constructing Lyapunov Functions Using Supervised Information}
Constructing Lyapunov functions has long been a significant challenge in dynamical systems and control, especially for applications involving complex and real-world tasks. This section addresses the challenge by leveraging supervised information to construct Lyapunov functions for applying to real-world tasks regardless of the network architecture or learning objective.

A critical step in constructing a Lyapunov function  $V$ is identifying the appropriate optimal state $h^*$. After determining $h^*$, the Lyapunov function $V$ is defined as
\begin{equation}
	{V_y}\left( h \right) := \frac{1}{2} \left\| {h - {h^*}} \right\|_2^2, \label{eq:Lya_fuction}
\end{equation}
which can be used to guarantee convergence of the dynamics to $h^*$. Recalling the training optimization problem for Neural ODEs, the optimal state $h^*$ needs to satisfy that the supervised loss is zero. For a classification task as an example, with supervised loss $\mathcal{L}_{cls}$ and a given input-output pair $(x, y)$, the optimal state $h^*$ can be obtained by solving a minimization problem as
 \begin{equation}
 	h^* :=\mathop{\arg \min}_{h \in H} \mathcal{L}_{cls}(\psi(h;\theta_{\psi}), y), \label{eq:h}
 \end{equation}
where $\psi$ is the output layer as defined in Equation \eqref{eq:2.3}. In general, the continuous differentiability of the output layer $\psi$ is easy satisfied, which ensures the feasibility of this approach.
 
Next we apply a multi-step gradient projection method for solving the optimization problem  \eqref{eq:h}, in which the choice of the initial point is crucial due to the fact that the unsuitable initial point may converge to a worthless $h^*$ for the dynamical system. By choosing a suitable initial point, we can increase the match between $ h^*$ and the dynamical system. Specifically, as shown in Figure \ref{fig:fullwidth-top}, we first use the solution trajectory $h(\cdot)$ in Equations \eqref{eq:2.1}-\eqref{eq:2.3} to determine the initial point $h(T)$, and then locate $h^*$ by the multi-step gradient projection method.

\subsection{ Fixed-Time Stabilizing Lyapunov Loss} \label{sec:3.2}
\begin{algorithm*}[ht]
	\caption{Robust learning Algorithms Based on Simulated Perturbation Sampling}
	\begin{algorithmic}[1]
		\State \textbf{Input:} Initial parameters $\theta$, learning rate $\eta_{1}$, number of iterations $N_{1}$, inner learning rate $\eta_{2}$, number of inner iterations $N_{2}$, number of samples $n_{\delta }$, $\varsigma_{max}$ denotes maximum disturbance radius, extracting feature dimensions $d$, time discretization resolution $\Gamma$, and  integral trajectory sampling spaced $t_0, t_1, \ldots, t_\Gamma$.
		\For{iteration $1$ to $N_{1}$}
		\State $(x, y) \sim D$  \Comment{Sample training data}
		\State $x_c=\phi \left( {x;{\theta _\phi }} \right)$  \Comment{Compute the extracted features of the input $x$}
		\State $h(t_i) \leftarrow \int_{t_{i-1}}^{t_i} f(h(\tau), x_c, \tau) \, d\tau + h(t_{i-1};\theta_{f}) \quad \forall \, 1 \le i \le \Gamma$ \Comment{Compute the integral trajectory in the input space}
		\For{iteration $1$ to $N_{2}$}
		\State $h^*=h(t_{\Gamma})$ \Comment{Assign initial value as $h(t_{\Gamma})$}
		\State ${h^*} \leftarrow {h^*} - \frac{{{\eta _2}{{\left| {h\left( {{t_\Gamma }} \right)} \right|}_2}}}{{{N_2}{{\left| {{\nabla _{{h^*}}}L(\psi ({h^*};{\theta _\psi }),y)} \right|}_2}}}{\nabla _{{h^*}}}L(\psi ({h^*};{\theta _\psi }),y)$ \Comment{Calculate the gradient and update $h^*$}
		\EndFor
		\State $h^{\delta_j}(t_i) \leftarrow \int_{t_{i-1}}^{t_i} f \left( h(\tau), x_c + \varsigma_j \|x_c\|_2 \frac{\delta_j}{\|\delta_j\|_2}, \tau \right) d\tau + h(t_{i-1}; \theta_f) \quad \delta_j \sim \mathcal{N}(0,1)^d, \varsigma_j \sim U(0,\varsigma_{\max}), \forall 1 \le j \le n_\delta, \forall 1 \le i \le \Gamma$ \Comment{Calculate perturbation integral trajectories and sampling}
		\State $\theta \leftarrow \theta - \eta_{1}(\nabla_\theta \sum_j \sum_i {\cal V}(x, y,h^{\delta_j}(t_i) , t_i)+{\nabla _{{\theta}}}{\cal L}_{cls}(\psi ({h(t_{\Gamma})};{\theta _\psi }),y))$ \Comment{Calculate total loss and compute gradient to update $\theta$} 
		\EndFor
		\State \textbf{Output:} Trained parameters $\theta$.
	\end{algorithmic}
\end{algorithm*}

In this subsection, we introduce a fixed-time stabilizing loss function for Neural ODEs, derived from the Lyapunov function  \( V \)  defined in Section \ref{sec:3.1}. The goal is to impose fixed-time stability on the learning dynamics by utilizing the stability condition in Equation \eqref{eq:FxS-stability}, leading to the formulation of FxTS-Loss. To begin, we define the pointwise version of this loss.
\begin{definition}[Point-wise FxTS-Loss] \label{def1}
	For a single input-output pair \( (x, y) \) and the Lyapunov function \( V_{y} : H \rightarrow \mathbb{R}_{\geq 0} \) defined in Equation \eqref{eq:Lya_fuction}, a point-wise FxTS-Loss defines as 
	\begin{align}\label{eq:PW_FxSL_loss}
		\mathcal{V}(x, y, h, t) := &\max \left\{ 0, \frac{\partial V}{\partial h}f(h, x_c, t;\theta_f) + \right. \nonumber \\
		& \left. \alpha_1 V_y(h)^{\gamma_1} + \alpha_2 V_y(h)^{\gamma_2}\right\},
	\end{align}
	where  $h \in H$, \(\alpha_1, \alpha_2 > 0\), \(\gamma_1 = 1 + \frac{1}{\mu}\) and \(\gamma_2 = 1 - \frac{1}{\mu}\) with \(\mu > 1\).
\end{definition}

This point-wise definition quantifies the degree of violation of the FxTS contraction condition at each local state \( h \), which can be interpreted as a local contraction property. To utilize Lemma \ref{lem:2.1}, we need to check that  \( \mathcal{V} \) is zero for all $H$ and ${t \in [0,1]}$ such that the dynamics converge to the loss-minimizing prediction within a user-defined fixed time.

Next, we define FxTS-Loss over the coupled distribution $\left\{ {(h\left( t \right),t)\left| {t \in [0,1]} \right.} \right\}$ (rather than the entire $H$ and ${t \in [0,1]}$) by integrating the pointwise loss over time.

\begin{definition}[FxTS-Loss] \label{def:FxSL_Loss}
	For the solution $h(\cdot)$ of dynamics defined in Equations \eqref{eq:2.1}-\eqref{eq:2.3}, and \( \mathcal{V} \) from Equation \eqref{eq:PW_FxSL_loss}, the FxSL-Loss is
	\begin{equation}\label{eq:in_loss}
		\mathcal{L}_{lyp}(\theta) := \mathbb{E}_{(x,y) \sim D} \left[ \int_{0}^{1} \mathcal{V}(x, y, h(t), t) \, dt \right].
	\end{equation}
	where $D$ is a  a dataset of input-output pairs.
\end{definition}
 
 With the help of Definition \ref{def:FxSL_Loss}, we end this subsection by giving our first main result. Theorem \ref{th:1}, shown below, demonstrates that minimizing FxTS-Loss guarantees fixed-time stability in the learned Neural ODE, with the proof given in \ref{sec:A1}.

\begin{theorem}\label{th:1}
	Under the setting of Definition \ref{def:FxSL_Loss}, if there exists a parameter \( \theta^* \in \Theta \) of the dynamical system that attains \( \mathcal{L}_{lyp}(\theta^*) = 0 \), then for almost everywhere  \( (x, y) \sim D \):
	\begin{enumerate}
		\item For any $t \in [0,1]$, the inference dynamical system satisfies $\mathcal{V}(x, y, h(t), t)=0 $.
		\item  All solution trajectories $h\left( t\right) $ of the dynamical system determined by Equations \eqref{eq:2.1}-\eqref{eq:2.3} satisfy
		\begin{equation}
		 \sup \{ T \geq 0 : h(t) = h^* \text{ for all } t \geq T \} \le \frac{\mu \pi}{2\sqrt{\alpha_1 \alpha_2}},
		\end{equation}
		where  $h^*$ is defined in Equation \eqref{eq:h} .
	\end{enumerate}
\end{theorem}

\subsection{ Robust learning Algorithms} \label{sec:learning_algorithms}

It is crucial to choose the appropriate set \( H \) and ensure that \( \mathcal{V} \) is equal to zero throughout ${\rm{H }} \times \left[ {0,1} \right]$. Although Theorem \ref{th:1} shows that the dynamical system convergences to the state \( h^* \) in a fixed time by minimizing the FxTS-Loss which formulates on the set of trajectories $\left\{ {(h\left( t \right),t)\left| {t \in [0,1]} \right.} \right\}$, it should be noticed that this stability depends on the integration trajectory, leading to a lack of robustness against input perturbations. Specifically, the choice of $H$ in Theorem  \ref{th:1} is so tight that a small input perturbation may change the integral trajectory of the dynamical system beyond the range of $H$. 

Therefore, such a challenge in designing robust learning algorithms lies in computing the inner integral in Equation \eqref{eq:in_loss}, which requires addressing two key problems:
\begin{itemize}
	\item Selecting an appropriate \( H \) that enables the learned Neural ODEs to exhibit robustness.
	\item Collecting sample points to approximate the FxTS-Loss over \( {\rm{H }} \times [0,1] \) with limited computational resources.
\end{itemize}

To tackle these key problems, we propose a simulated perturbation sampling method. As illustrated in the state space visualization in Figure \ref{fig:fullwidth-top}, we simulate input perturbations by perturbation extracting perturbation features to generate the post-perturbation integral trajectory (as shown by the orange curve in the figure) and subsequently sample this trajectory. The pseudo-code for the robust learning algorithm is provided in Algorithm \ref{as:1}. This method effectively collects sampling points in the critical regions of the Neural ODE, thereby enhancing both robustness and accuracy.

\subsection{ Adversarial Robustness}
\begin{figure}[t]
	\centering
	\includegraphics[width=0.9\columnwidth]{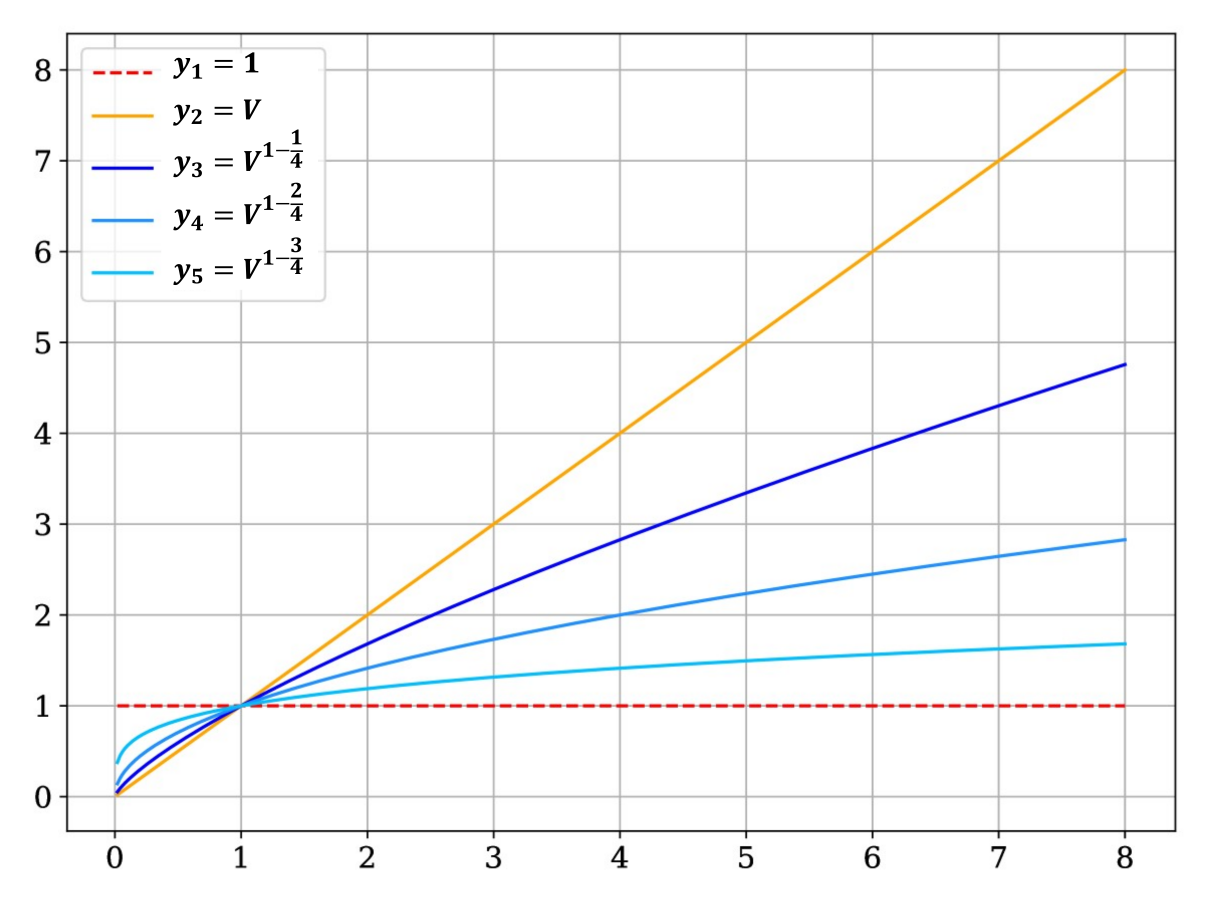}
	\vspace{-15pt}
	\caption{Plots of Scale Functions (setting $a_1,a_2,\delta=1$, and $\mu=4$). The red dotted line indicates error-free, i.e., multiplied by $1$ on $\delta$. The orange solid line indicates the scaling technique in \cite{garg2021advances}, which multiplying $V$ on $\delta$. The remaining solid lines indicate ours.}
	\vspace{-10pt}
	\label{fig:scale-functions}
\end{figure}

It is well known that small perturbations at the input of an unstable dynamical system will lead to large distortions in the system's output. This section will explore the provable adversarial robustness obtained by minimizing the FxTS-Loss, which utilizes robust control based on Lyapunov stability theory. Robust control fundamentally ensures that a system remains stable when perturbed, and the notion is closely related to the study of adversarial robustness in machine learning \cite{WANG2024106087,LI2023177}. 

Firstly, to obtain adversarial robustness guarantees, we develop more precise time upper bound estimation for non-vanishingly perturbed systems \cite[Section 9.2]{Khalil:1173048}.  Specifically, considering the dynamical system in Equations \eqref{eq:2.1}-\eqref{eq:2.3}, for an input perturbation, we can utilize the Lipschitz continuity from Assumption \ref{as:1} to transform the perturbed dynamics into the original dynamics with a bounded nonvanishing perturbation. 

Moreover, we relax FxTS Lyapunov conditions in  Equation \eqref{eq:FxS-stability} by allowing a constant term to appear in the upper bound of the time derivative of the Lyapunov function. The FxTS Lyapunov condition is considered for a positive definite, proper, continuously differentiable \( V : \mathbb{R}^n \rightarrow \mathbb{R} \) as follows
\begin{equation} \label{eq:A8}
	\dot{V}(x) \leq -\alpha_1 V(x)^{\gamma_1} - \alpha_2 V(x)^{\gamma_2} + \delta,
\end{equation}
where \(\alpha_1, \alpha_2 > 0\), \(\delta > 0 \), \(\gamma_1 = 1 + \frac{1}{\mu}\), \(\gamma_2 = 1 - \frac{1}{\mu}\) and \(\mu > 1\). Remark that if \(\delta=0\), then Lyapunov condition \eqref{eq:A8}  is reduced  to Equation \eqref{eq:FxS-stability}.  Proposition \ref{propos:1}, shown below, provides the expression for the time upper bound estimation of the equation \eqref{eq:A8}, with the proof given in  \ref{sec:A2}.
\begin{proposition} \label{propos:1}
	Let \( V_0, \alpha_1, \alpha_2, \delta, v > 0 \), and $\gamma >1$ such that $\alpha_{1}v^{\mu+1}+\alpha_2v^{\mu-1}=\delta$. Let \( \gamma_1 = 1 + \frac{1}{\mu} \) and \( \gamma_2 = 1 - \frac{1}{\mu} \) with \( \mu > 1 \). Define
	\begin{equation}
		I: = \int_{{V_0}}^{\bar V} {\frac{{dV}}{{ - {\alpha _1}{V^{{\gamma _1}}} - {\alpha _2}{V^{{\gamma _2}}} + \delta }}}. 
	\end{equation}
	Then, for all $V_0 \ge \bar{V}=(\gamma v)^{\mu}$, the following conclusions hold:
	\begin{enumerate}
		\item  For $1 < \mu < 2$, one has 
		\begin{equation}
			I \le \frac{\mu }{{2{\alpha _1}v}}\left[ {\ln\left( {\left| {\frac{{1 + \gamma }}{{\gamma  - 1}}} \right|} \right) + \ln\left( {\left| {\frac{{V_0^{\frac{1}{\mu }} - v}}{{V_0^{\frac{1}{\mu }} + v}}} \right|} \right)} \right]. \nonumber
		\end{equation}
		
		\item  For $2 \le \mu < 3$, one obtains
		\begin{align}
			I \leq & k_1
			\left[ \ln {\left| {\frac{{\left( {{\gamma ^2} + \gamma  + 1} \right){v^2} + \frac{{{\alpha_2}}}{{{\alpha_1}}}}}{{{{\left( {\left( {\gamma  - 1} \right)v} \right)}^2}}}} \right|} \right. \nonumber \\
			& + \ln {\left| {\frac{{{{\left( {{V_0}^{\frac{1}{\mu }} - v} \right)}^2}}}{{{V_0}^{\frac{2}{\mu }} + v{V_0}^{\frac{1}{\mu }} + {v^2} + \frac{{{\alpha_2}}}{{{\alpha_1}}}}}} \right|}  \nonumber \\
			& \left. + k_2
			\left( \arctan {\frac{{2{V_0}^{\frac{1}{\mu }} + v}}{{\sqrt {3{v^2} + \frac{{4{\alpha_2}}}{{{\alpha_1}}}} }}} - \arctan {\frac{{\left( {2\gamma  + 1} \right)v}}{{\sqrt {3{v^2} + \frac{{4{\alpha_2}}}{{{\alpha_1}}}} }}} \right)  \right] ,\nonumber
		\end{align}
		where $k_1= \frac{\mu v}{6 \alpha_1 v^2 + 2 \alpha_2} $, $ k_2=\frac{6v + \frac{4 \alpha_2}{\alpha_1 v}}{\sqrt{3v^2 + \frac{4 \alpha_2}{\alpha_1}}}$.

		\item  For $\mu \ge 3$, it holds that
			\begin{align}
				I \leq &\ k_3 \left[ \frac{v}{2} \left( \ln\left| \frac{1 + \gamma}{1 - \gamma} \right| + \ln\left| \frac{V_0^{\frac{1}{\mu}} - v}{V_0^{\frac{1}{\mu}} + v} \right| \right) \right. \nonumber \\
				& \left. + k_4 \left( \arctan \left( \frac{V_0^{\frac{1}{\mu}}}{k_4} \right) - \arctan \left( \frac{\gamma v}{k_4} \right) \right) \right], \nonumber
			\end{align}

		where $k_3=\frac{\mu }{{2{\alpha _1}{v^2} + {\alpha_2}}}$, $k_4=\sqrt {{v^2} + \frac{{{\alpha _2}}}{{{\alpha _1}}}}$.
	\end{enumerate}
\end{proposition}

It is worth mentioning that Proposition \ref{propos:1} provides more precise estimations for both the upper bound on time $I$ and the boundary $\bar{V}$ than the ones in \cite{garg2021advances}. In fact, \citet{garg2021advances} scaled Equation \eqref{eq:A8} by multiplying $V$ on $\delta$ to solve the failure in integration due to the introduction of a constant term. This means that the convergence boundary also needs to satisfy $\bar{V}>1$, besides multiplying $V$ is not the best choice for scaling. In Proposition \ref{propos:1}, we multiply ${V^{1 - \frac{n}{\mu }}}$ for more accurate estimation and by introducing ${v^{1-\mu}}$ to eliminate the condition $V>1$.  In Figure \ref{fig:scale-functions}, we visualize  the different scaling techniques and can see that our scaling technique is having a tighter effect. Notably, when $\mu$ is selected as $2,3$, our estimation results are precise and error-free. 

\begin{table*}[ht]
	\centering
	\caption{ Classification error and robustness against stochastic noises. ``Original" denotes the classification error on original data. ``Avg" denotes the average error of the five random noise-perturbed images and the originals over the three datasets. The best result is in bold. Unit: \%.} 
	\label{table:acc_noise}
	\begin{tabular}{ccccccccc}
		\toprule
		Dataset                   & Method     & Orignal        & Gaussian       & Shot           & Impulse        & Speckle          & Motion    & Avg    \\ \hline
		\multirow{4}{*}{Fashion-MNIST} & ResNet18   & 8.44            & 16.37           & 14.47           & 41.15           & 9.18            & 12.35  & 16.99  \\
		& Neural ODE & 8.33            & 15.96           & 14.90           & 40.69           & 9.50            & 12.18        & 16.93  \\
		& LyaNet     & 8.12            & 15.09           & \textbf{13.22}  & 40.43           & \textbf{9.08}   & 12.83        & 16.46   \\
		& FxTS-Net   & \textbf{8.08}   & \textbf{14.68}  & 14.10           & \textbf{32.60}  & 9.15     & \textbf{11.74}    & \textbf{15.06} \\ \hline
		\multirow{4}{*}{CIFAR10}  & ResNet18        & 18.35            & 51.38             & 55.93          & 60.60            & 24.88 & 26.87   & 39.67  \\
		& Neural ODE      & 17.64            & 48.42             & 53.04          & 57.32            & 23.94      & 24.92         & 37.55   \\
		& LyaNet          & 17.53            & 46.78             & 50.40          & 54.84            & 23.35    & 24.65  & 36.26  \\
		& FxTS-Net        & \textbf{16.21}   & \textbf{43.96}    & \textbf{48.90} & \textbf{52.64}   & \textbf{21.63}  & \textbf{24.47}       & \textbf{34.64}    \\ \hline
		\multirow{4}{*}{CIFAR100} & ResNet18        & 41.43            & 71.84             & 73.98          & 79.75            & 45.75  & 46.96   & 59.95        \\
		& Neural ODE      & 41.27            & 70.94             & 73.03          & 79.22            & 45.61  & 46.27    & 59.39       \\
		& LyaNet          & 40.48            & 69.41             & 71.90          & 77.25            & 45.68   & 46.30     & 58.50      \\
		& FxTS-Net        & \textbf{39.07}   & \textbf{67.92}    & \textbf{69.84} & \textbf{76.47}   & \textbf{43.25}  & \textbf{44.42}   & \textbf{56.83} \\  \bottomrule
	\end{tabular}
\end{table*}

Secondly, we use Proposition \ref{propos:1} to further show that optimizing FxTS-Loss can obtain robustness guarantees. Let $(x,y)$ denote the input-output pair, and $h_{x}(\cdot)$ be the solution of the dynamical system defined in Equations \eqref{eq:2.1}-\eqref{eq:2.3}. Consider $\tilde{x}$ to be a perturbed version of the input data such that
\begin{equation}
	\left\| {x - \tilde x} \right\|_2^2\ \leq \rho, \label{eq:16}
\end{equation}
where $\rho > 0$ represents the degree of perturbation. Suppose that we consider the supervised learning task mentioned in Section 2 and use $ h_{\tilde{x}}(\cdot)$ to represent the perturbed version of $ h_{x}(T)$.

\begin{theorem}\label{the:2}
	Let $(x, y) $ be an input data pair and $\tilde{x}$ be the respective perturbed version, $x$ and $\tilde{x}$ satisfy Equation \eqref{eq:16}, $L_V, L_f, L_{\phi}, L_{\psi}$ be the Lipschitz constants from Assumption \ref{as:1}, $\gamma >1$, $\tilde{y}(t)=\psi \left( {h_{\tilde{x}}\left( t \right);{\theta _\psi }} \right)$, and $v$ satisfy ${\alpha _1}{v^{u + 1}} + {\alpha _2}{v^{\mu  - 1}} = L\rho$ with $L={{L_V}{L_f}{L_{\phi}}}$. Then  $\left| {\tilde y(t) - y} \right| \le {\left( {\gamma v} \right)^\mu }{L_\psi }$ for all $ t \ge T$, where $T$ has the following formation:
	\begin{enumerate}
		\item  For $1 \le \mu \le 2$, we have
		\begin{equation}
			T = \frac{\mu }{{2{\alpha _1}v}}\ln\left( {\left| {\frac{{(1 + \gamma )v}}{{(\gamma  - 1)v}}} \right|} \right). \nonumber
		\end{equation}
		\item  For $2 \le \mu \le 3$, we get

		\begin{align}
			T = & \frac{\mu v}{6\alpha_1 v^2 + 2\alpha_2} \left[ \ln \left| \frac{(\gamma^2 + \gamma + 1)v^2 + \frac{\alpha_2}{\alpha_1}}{((\gamma - 1)v)^2} \right| \right. \nonumber \\
			& \left. + \frac{6v + \frac{4\alpha_2}{\alpha_1 v}}{\sqrt{3v^2 + \frac{4\alpha_2}{\alpha_1}}} \left( \frac{\pi}{2} - \arctan \left( \frac{(2\gamma + 1)v}{\sqrt{3v^2 + \frac{4\alpha_2}{\alpha_1}}} \right) \right) \right]. \nonumber
		\end{align}
		\item  For $\mu \ge 3$, one has
		\begin{align}
			T = & \frac{\mu}{2\alpha_1 v^2 + \alpha_2} \left[ \frac{v}{2} \left( \ln \left| \frac{1 + \gamma}{1 - \gamma} \right| \right) \right. \nonumber \\
			& \left. + \sqrt{v^2 + \frac{\alpha_2}{\alpha_1}} \left( \frac{\pi}{2} - \arctan \left( \frac{\gamma v}{\sqrt{v^2 + \frac{\alpha_2}{\alpha_1}}} \right) \right) \right]. \nonumber
		\end{align}
	\end{enumerate}

\end{theorem}

See \ref{sec:A3} for the proof. Basically, Theorem \ref{the:2} provides the stabilization to a neighborhood of origin (the radius of the neighborhood is proportional to the perturbation magnitude) at a fixed time for the difference between the perturbed outputs and label $y$. It is worth saying that the radius of the neighborhood can be effectively reduced by setting the appropriate $\alpha_{1}$ and $\alpha_{2}$ in the learning.  In contrast, the lack of such stability in other Neural ODEs can yield dramatically different solutions from small changes in the input data, causing dramatic performance degradation.

\section{Experiments}
In this section, we conduct a series of comprehensive experiments to evaluate the performance of our proposed FxTS-Net model. Our analysis aims to address the following key research questions:
\begin{itemize}
	\item RQ1: How does FxTS-Net perform in terms of classification accuracy compared to baseline models?
	\item RQ2: How robust is FxTS-Net against various types of perturbations, including stochastic noise and adversarial attacks?
	\item RQ3: What impact does the fixed-time stability guarantee have on the decision boundaries of FxTS-Net compared to traditional Neural ODEs?
	\item RQ4: How do different parameters affect the performance of FxTS-Net?
\end{itemize}

Through these experiments, we provide a comprehensive analysis of the model's performance and robustness across various settings, highlighting the advantages of incorporating fixed-time stability.

\subsection{Experimental settings}

This subsection provides a detailed overview of the experimental setup, including the datasets used, training configurations, and evaluation metrics employed to assess the performance and robustness of FxTS-Net.

\begin{table*}[ht]
	\centering
	\caption{Robustness against adversarial attacks with different attack radii. The best results are shown in bold. Unit: \%.} 
	\label{table:attac}
	\begin{tabular}{cccccccccc}
		\toprule
		Dataset                   & Method               & \multicolumn{2}{c}{FGSM}        & \multicolumn{2}{c}{BIM}         & \multicolumn{2}{c}{PGD}         & \multicolumn{2}{c}{APGD}\\\cmidrule(r){3-4} \cmidrule(r){5-6} \cmidrule(r){7-8} \cmidrule(r){9-10}
		\multicolumn{1}{l}{}      & \multicolumn{1}{l}{} & 8/255          & 16/255         & 8/255          & 16/255         & 8/255          & 16/255         & 8/255          & 16/255    \\ \hline
		\multirow{4}{*}{Fashion-MNIST} & ResNet18             & 52.46           & 62.84           & 23.84           & 27.46           & 36.52           & 37.44           & 39.56   & 39.81    \\
		& Neural ODE           & 42.40           & 51.20           & 21.85           & 23.93           & 30.75           & 31.67           & 32.91           & 33.01           \\
		& LyaNet               & 30.76           & 33.77           & 11.97           & 12.54           & 27.87           & 28.45           & 29.74           & 30.27           \\
		& FxTS-Net             & \textbf{25.85}  & \textbf{26.06}  & \textbf{8.28}   & \textbf{8.41}   & \textbf{24.72}  & \textbf{27.14}  & \textbf{26.94}  & \textbf{27.21}  \\ \hline
		\multirow{4}{*}{CIFAR10}  & ResNet18             & 49.39          & 60.45          & 35.28          & 39.13          & 37.58          & 38.27          & 40.93          & 42.17          \\
		& Neural ODE           & 49.86          & 59.90          & 37.09          & 40.94          & 38.20          & 39.96          & 42.01          & 44.05          \\
		& LyaNet               & 44.23          & 51.03          & 27.11          & 28.12          & 38.30          & 39.28          & 41.58          & 42.63          \\
		& FxTS-Net             & \textbf{37.37} & \textbf{38.33} & \textbf{19.15} & \textbf{19.96} & \textbf{36.03} & \textbf{37.25} & \textbf{39.48} & \textbf{40.17} \\ \hline
		\multirow{4}{*}{CIFAR100} & ResNet18             & 84.73          & 87.20          & 62.72          & 65.52          & 82.34          & 83.31          & 84.98          & 85.13         \\
		& Neural ODE           & 84.69          & 87.03          & 62.86          & 65.34          & 81.96          & 82.34          & 84.65          & 84.97          \\
		& LyaNet               & 81.24          & 81.84          & \textbf{55.92} & \textbf{56.33} & 80.45          & 80.97         & 83.38          & 83.93          \\
		& FxTS-Net             & \textbf{80.00} & \textbf{80.72} & 56.77          & 57.54          & \textbf{78.82} & \textbf{79.72} & \textbf{81.49} & \textbf{81.58} \\  \bottomrule
	\end{tabular}
\end{table*}

\noindent \textbf{Datasets.}\quad We evaluate our method on three classification datasets: Fashion-MNIST \cite{xiao2017fashion}, CIFAR-10, and CIFAR-100 \cite{krizhevsky2009learning}. Fashion-MNIST consists of 70,000 grayscale images (28×28 pixels) across 10 classes, with 60,000 used for training and 10,000 for testing. CIFAR-10 and CIFAR-100 contain 60,000 color images (32×32 pixels), with 50,000 images for training and 10,000 for testing, categorized into 10 and 100 classes, respectively.

\noindent  \textbf{Training Details.}\quad  Our implementation is based on the open-source codes.\footnote{https://github.com/ivandariojr/LyapunovLearning} The Runge-Kutta method \cite{NEURIPS2018_69386f6b} is still used as numerical solvers in the process of forward propagation. Following \cite{pmlr-v162-rodriguez22a} to simplify tuning, we trained models using Nero \cite{pmlr-v139-liu21c} with a learning rate of 0.01, a batch size 64 for LyaNet and Ours, and 128 for other models. All models were trained for a total of 120 epochs. In our model, we set inner learning rate $\eta_{2}=2$, inner iterations $N_{2}=3$, samples $n_{\delta }=102$, and  time discretization resolution $\Gamma=5$.

\noindent \textbf{Evaluation Details.}\quad Model performance is measured by classification error on the test sets of Fashion-MNIST, CIFAR-10, and CIFAR-100. Following \cite{CUI2023576, YAN2020On}, robustness is assessed using two types of perturbations: stochastic noise and adversarial attacks. For stochastic noise, we apply Gaussian, shot, impulse, speckle, and motion noise \cite{hendrycks2018benchmarking}. For adversarial robustness, we consider the Fast Gradient Sign Method (FGSM) \cite{goodfellow2014explaining}, Basic Iterative Method (BIM) \cite{kurakin2018adversarial}, Projected Gradient Descent (PGD) \cite{madry2018towards}, and Auto Projected Gradient Descent (APGD) \cite{pmlr-v119-croce20b}.

We compare three models: 
\begin{itemize} 
	\item \textbf{ResNet-18} \cite{He_2016_CVPR}: A conventional convolutional neural network where the residual connections approximate Eulerian integrated dynamical systems. 
	\item \textbf{Neural ODEs} \cite{pmlr-v162-rodriguez22a}: A data-controlled Neural ODE framework using ResNet-18 as a feature extractor. 
	\item \textbf{LyaNet} \cite{pmlr-v162-rodriguez22a}: A stable Neural ODE model that guarantees exponential stability via the Lyapunov exponential stability condition. \end{itemize}

\subsection{Accuracy and Robustness Against Stochastic Noises}

This subsection focuses on the classification errors on standard and noise-perturbed images, with a fixed maximum evolution time of $1$. The results, detailed in Table \ref{table:acc_noise}, show that FxTS-Net consistently outperforms the benchmark models across all datasets. Specifically, FxTS-Net achieves a 0.04\%, 1.32\%, and 1.41\% improvement in standard classification accuracy than the best one of other models on Fashion-MNIST, CIFAR-10, and CIFAR-100, respectively. This demonstrates that our fixed-time stable learning framework enables Neural ODEs to make accurate predictions within a pre-defined time, effectively enhancing classification performance.

Moreover, FxTS-Net exhibits improved robustness on noise-perturbed data, achieving the best average performance across all random noise types. Although it lags slightly behind LyaNet by 0.07\% in speckle noise and 0.88\% in shot noise for Fashion-MNIST, FxTS-Net remains highly competitive. These results emphasize the robustness of FxTS-Net to small input perturbations, validating our model's ability to regulate the dynamic behavior of Neural ODEs and enhance stability. This aligns with our initial motivation to improve the robustness of Neural ODEs against random noise perturbations.

\subsection{Robustness Against Adversarial Attacks}

Similar to our evaluation of randomly noisy images, we evaluated the adversarial robustness of FxTS-Net and the benchmark models within a fixed evolutionary time Each model was tested on Fashion-MNIST, CIFAR-10, and CIFAR-100 datasets, using adversarial perturbations with attack radii of $\epsilon = 8/255$ and $\epsilon = 16/255$ from FGSM, BIM, PGD, and APGD attacks, as detailed in Table \ref{table:attac}.

Specifically, on Fashion-MNIST and CIFAR-10, FxTS-Net consistently outperforms the other models across all attacks. For CIFAR-100, FxTS-Net demonstrates superior robustness across FGSM, PGD, and APGD attacks and only slightly lags behind LyaNet by 0.85\% and 1.21\% in the BIM attack. The observed improvement underscores the efficacy of our robust learning algorithm, which effectively captures key sample points in critical regions to compute the FxTS-Loss. This facilitates the local-global structure alignment of the Lyapunov condition, confirming the importance of stability guarantees within pre-defined fixed times for Neural ODEs. The results also support our theoretical analysis, showing that efficiently approximating FxTS-Loss enhances robustness to adversarial attacks.

\begin{figure*}[h] 
	\centering
	\includegraphics[width=\textwidth]{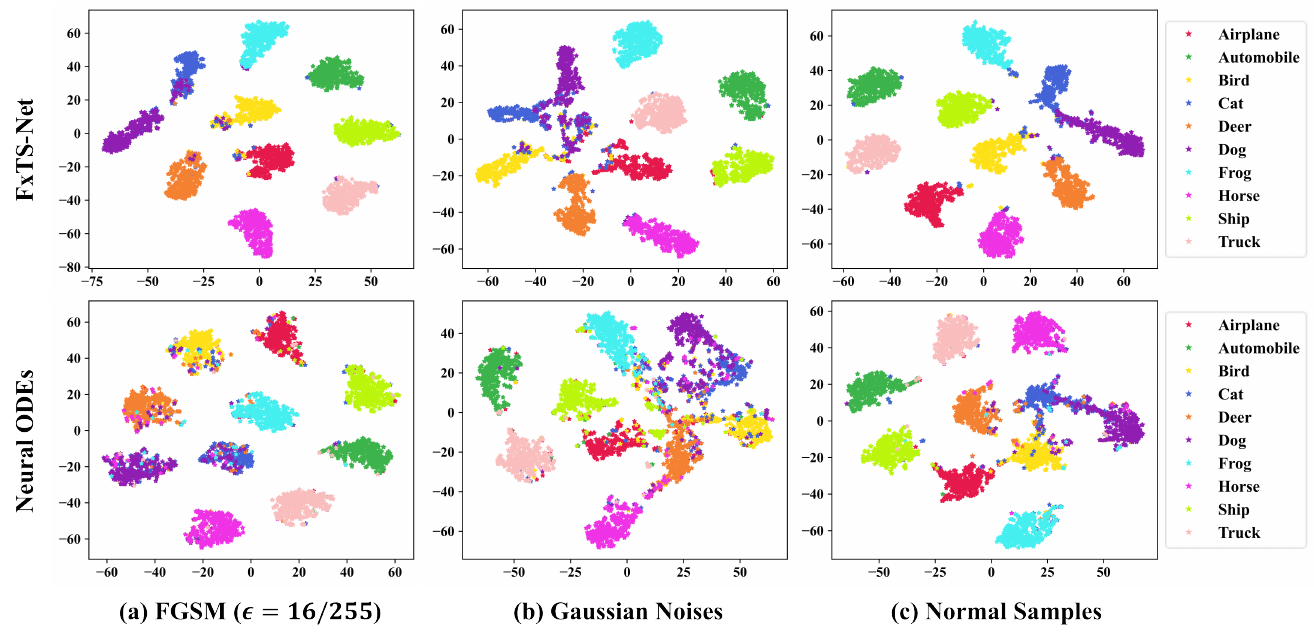} 
	\caption{TSNE visualisation results on classification features of the Neural ODE and FxTS-Net using CIFAR10 with 400 random images per category. From left to right, FGSM attack perturbation, Gaussian noise perturbation and normal data.} 
	\label{fig:tsne}
\end{figure*}

\begin{table*}[th]
	\centering
	\caption{ Classification error and robustness to perturbations for the ablation study of the parameter $\varsigma_{max}$, where the adversarial attack perturbations $\epsilon=8/255$. ``Original" denotes the classification error on original data. The best result is in bold. ``Avg" denotes the average error. Unit: \%.} \label{table:ab2}
	\begin{tabular}{cccccccccc}
		\toprule
		Parameters            & Orignal        & Gaussian       & Impulse        & Speckle        & Motion         & FGSM           & BIM            & APGD         & Avg            \\ \hline
		$\varsigma_{max}=0.4$ &  \textbf{16.06}          & 46.44          & 55.39          & 22.80          & 25.16          & 37.90          & 19.65          & 39.78          & 32.90          \\
		$\varsigma_{max}=0.8$ & 16.15          & 44.82          & 53.27          & 22.03          & 25.45  & 37.67          & 19.58          & 39.92          & 32.36          \\
		$\varsigma_{max}=1.2$ & 16.21 & \textbf{43.96} & \textbf{52.64} & \textbf{21.63} & 24.47          & \textbf{37.37} & \textbf{19.15} & \textbf{39.48} & \textbf{31.86} \\
		$\varsigma_{max}=1.6$ & 16.74          & 44.28          & 52.74          & 21.88          &  \textbf{23.81}         & 37.53          & 19.42          & 39.75          & 32.01      \\  \bottomrule    
	\end{tabular}
\end{table*}

\subsection{Inspecting Decision Boundaries}

In this subsection, we further investigate the impact of introducing fixed-time stability guarantees on the decision boundaries learned by Neural ODEs. We perform a TSNE clustering analysis on the classification features of Neural ODEs and FxTS-Net using adversarial perturbations from FGSM attacks, Gaussian noise, and clean CIFAR-10 data. The results, shown in Figure \ref{fig:tsne}, highlight clear differences between the models. For standard data, FxTS-Net demonstrates a more distinct separation between clusters of different categories compared to Neural ODEs. Notably, under Gaussian noise, the decision boundaries of Neural ODEs become highly blurred, whereas FxTS-Net maintains sharp boundaries. This illustrates how the fixed-time stability framework enhances the model's decision boundaries and mitigates boundary ambiguity under noisy conditions.

Moreover, while Gaussian noise causes boundary-blurring, FGSM perturbations do not blur the decision boundaries but directly confuse the classification outcomes. Even under adversarial attacks, FxTS-Net consistently delivers better classification performance than Neural ODEs, which strongly supports the effectiveness of our approach. In summary, FxTS-Net excels in maintaining stable decision boundaries and classification accuracy across various data perturbations, validating the role of fixed-time stability in improving both the robustness and classification performance of Neural ODEs.

\subsection{Ablation Study of Parameters}
\label{sec:ASP}
In this subsection, we conduct ablation experiments to explore the effects of different parameters in FxTS-Net, using CIFAR10 as the dataset. Firstly, we conduct an ablation study of the algorithm's maximum perturbation radius $\varsigma_{max}$. Secondly, we conduct experiments for the parameters $\mu, \alpha_{1}, \alpha_{2}$ in the FxTS-Loss function to assess their impact on the robustness under various conditions. Additional ablation study results are described in Appendix B. These experiments aim to reveal the key role of each parameter on the stability and noise robustness of FxTS-Net, to guide parameter selection.

The first part of the experiment investigates the effect of the perturbation radius $\varsigma_{max}$. The results in Table \ref{table:ab2} demonstrate that smaller values of $\varsigma_{max}$ (e.g., 0.4) can lead to better performance on clean data but at the expense of robustness under adversarial conditions. Conversely, increasing $\varsigma_{max}$ to 1.2 achieves the best overall performance across both original and adversarial data, with the highest average accuracy. This indicates that moderate perturbations during training encourage the model to generalize well across different types of inputs, thereby enhancing its robustness.  However, for larger values of $\varsigma_{max}$, FxTS-Net exhibits diminishing returns, indicating an optimal disturbance range for balancing accuracy and robustness.

\begin{table*}[th]
	\centering
	\caption{  Classification error and robustness to perturbations for the ablation study of parameters $\mu, \alpha_{1}, \alpha_{2}$, where the adversarial attack perturbations $\epsilon=8/255$. ``Original" denotes the classification error on original data. The best result is in bold. Unit: \%.} \label{table:ab1}
	\begin{tabular}{lccccccccc}
		\toprule
		\multicolumn{1}{c}{Parameters}              & Orignal        & Gaussian       & Impulse        & Speckle        & Motion         & FGSM           & BIM            & APGD         & Avg           \\ \hline
		$\mu=2.0, \alpha_{1}=2.0, \alpha_{2}=4.94$  & 17.10          & 44.87          & 52.71          & 22.56          & 25.51          & 37.86          & 19.72          & 39.79          & 32.52          \\
		$\mu=3.0, \alpha_{1}=2.0, \alpha_{2}=11.20$ & 16.57 & 45.24          & 53.60          & 22.09          & 24.66          & 37.59          & \textbf{19.13} & 39.87          & 32.34          \\
		$\mu=2.0, \alpha_{1}=6.0, \alpha_{2}=1.65$  & 16.94          & 45.98          & 54.53          & 22.48          & 25.29          & 37.86          & 19.84          & 40.03          & 32.87          \\
		$\mu=3.0, \alpha_{1}=6.0, \alpha_{2}=3.71$  & 17.41          & 44.88          & 53.47          & 22.35          & 25.47          & 38.03          & 20.02          & 39.74          & 32.67          \\
		$\mu=2.0, \alpha_{1}=10.0, \alpha_{2}=1.00$ & \textbf{16.21} & \textbf{43.96} & \textbf{52.64} & \textbf{21.63} & \textbf{24.47} & \textbf{37.37} & 19.15          & \textbf{39.48} & \textbf{31.86} \\
		$\mu=3.0, \alpha_{1}=10.0, \alpha_{2}=2.23$ & 16.70          & 44.11          & 52.78          & 22.18          & 25.11          & 37.56          & 19.86          & 39.59          & 32.24         \\  \bottomrule 
	\end{tabular}
\end{table*}

In the second part, we explore how different combinations of $\mu$, $\alpha_{1}$, and $\alpha_{2}$ affect model performance. As shown in Table \ref{table:ab1} , the parameter combination $\mu=2.0, \alpha_{1}=10.0, \alpha_{2}=1.00$ shows the best results on both raw and attack data (including FGSM and motion blur), especially on the average error of 31.86\%. It is also shown that the larger value of $\alpha_1$, compared to $\alpha_2$, leads to better accuracy and robustness of the model, while the variation of $\mu$ between 2.0 and 3.0 has less impact on the final results. 

In summary, the ablation study highlights the importance of tuning both the loss function parameters and disturbance radius to balance between accuracy and robustness. The results demonstrate that the combination of $\alpha_{1}=10.0$, $\alpha_{2}=1.00$, and $\varsigma_{max}=1.2$ provides the best overall performance, particularly in adversarial settings. 

\section{Conclusion}

In this work, we introduce FxTS-Net, a framework for ensuring fixed-time stability of Neural ODEs by incorporating novel FxTS Lyapunov losses. We develop a more precise time upper bound estimation for bounded non-decreasing perturbation systems, providing a theoretical foundation for extending Neural ODEs to fixed-time stabilization and addressing their theoretical guarantees regarding adversarial robustness.

Furthermore, we present a method for constructing Lyapunov functions by utilizing supervised information during training, which can be adapted to the requirements of different tasks and network structures. We also provide learning algorithms, including capturing sample points in the critical region to approximate the FxTS-Loss, which promotes the local-global structure of the Lyapunov condition. Our experimental results demonstrated the effectiveness of FxTS-Net in improving adversarial robustness while maintaining competitive classification accuracy.

While our approach demonstrates competitive classification accuracy and robustness, several avenues remain for future exploration. One is the extension of FxTS-Net to more complex and high-dimensional tasks, such as video classification or 3D spatial data modeling. The other is the exploration of FxTS-Net to  integrate with some known adversarial defense mechanisms, which ensures stronger robustness. 

\section*{Declaration of competing interest}
The authors declare that they have no known competing financial interests or personal relationships that could have appeared to influence the work reported in this paper.

\section*{Acknowledgement}
This work was supported by the National Natural Science Foundation of China (12171339, 12471296).


\bibliographystyle{elsarticle-num-names} 
\bibliography{cas-refs}

\begin{thebibliography}{34}
\providecommand{\natexlab}[1]{#1}
\providecommand{\url}[1]{\texttt{#1}}
\providecommand{\urlprefix}{URL }
\expandafter\ifx\csname urlstyle\endcsname\relax
  \providecommand{\doi}[1]{doi:\discretionary{}{}{}#1}\else
  \providecommand{\doi}[1]{doi:\discretionary{}{}{}\begingroup
  \urlstyle{rm}\url{#1}\endgroup}\fi
\providecommand{\bibinfo}[2]{#2}

\bibitem[{Haber and Ruthotto(2017)}]{Haber_2018}
\bibinfo{author}{E.~Haber}, \bibinfo{author}{L.~Ruthotto},
  \bibinfo{title}{Stable architectures for deep neural networks},
  \bibinfo{journal}{Inverse Problems}
  \bibinfo{volume}{34}~(\bibinfo{number}{1}) (\bibinfo{year}{2017})
  \bibinfo{pages}{014004}.

\bibitem[{Ruthotto and Haber(2020)}]{Ruthotto2020}
\bibinfo{author}{L.~Ruthotto}, \bibinfo{author}{E.~Haber}, \bibinfo{title}{Deep
  Neural Networks Motivated by Partial Differential Equations},
  \bibinfo{journal}{Journal of Mathematical Imaging and Vision}
  \bibinfo{volume}{62}~(\bibinfo{number}{3}) (\bibinfo{year}{2020})
  \bibinfo{pages}{352--364}.

\bibitem[{Kidger et~al.(2021{\natexlab{a}})Kidger, Chen, and
  Lyons}]{pmlr-v139-kidger21a}
\bibinfo{author}{P.~Kidger}, \bibinfo{author}{R.~T.~Q. Chen},
  \bibinfo{author}{T.~J. Lyons}, \bibinfo{title}{``Hey, that’s not an ODE":
  Faster ODE Adjoints via Seminorms}, in: \bibinfo{booktitle}{Proceedings of
  the 38th International Conference on Machine Learning}, vol.
  \bibinfo{volume}{139}, \bibinfo{pages}{5443--5452},
  \bibinfo{year}{2021}{\natexlab{a}}.

\bibitem[{Oh et~al.(2024)Oh, Lim, and Kim}]{oh2024stable}
\bibinfo{author}{Y.~Oh}, \bibinfo{author}{D.~Lim}, \bibinfo{author}{S.~Kim},
  \bibinfo{title}{Stable Neural Stochastic Differential Equations in Analyzing
  Irregular Time Series Data}, in: \bibinfo{booktitle}{The Twelfth
  International Conference on Learning Representations}, \bibinfo{year}{2024}.

\bibitem[{Kidger et~al.(2020)Kidger, Morrill, Foster, and
  Lyons}]{NEURIPS2020_4a5876b4}
\bibinfo{author}{P.~Kidger}, \bibinfo{author}{J.~Morrill},
  \bibinfo{author}{J.~Foster}, \bibinfo{author}{T.~Lyons},
  \bibinfo{title}{Neural Controlled Differential Equations for Irregular Time
  Series}, in: \bibinfo{booktitle}{Advances in Neural Information Processing
  Systems}, vol.~\bibinfo{volume}{33}, \bibinfo{pages}{6696--6707},
  \bibinfo{year}{2020}.

\bibitem[{Yildiz et~al.(2019)Yildiz, Heinonen, and
  Lahdesmaki}]{NEURIPS2019_99a40143}
\bibinfo{author}{C.~Yildiz}, \bibinfo{author}{M.~Heinonen},
  \bibinfo{author}{H.~Lahdesmaki}, \bibinfo{title}{ODE2VAE: Deep generative
  second order ODEs with Bayesian neural networks}, in:
  \bibinfo{booktitle}{Advances in Neural Information Processing Systems},
  vol.~\bibinfo{volume}{32}, \bibinfo{year}{2019}.

\bibitem[{Kidger et~al.(2021{\natexlab{b}})Kidger, Foster, Li, and
  Lyons}]{pmlr-v139-kidger21b}
\bibinfo{author}{P.~Kidger}, \bibinfo{author}{J.~Foster},
  \bibinfo{author}{X.~Li}, \bibinfo{author}{T.~J. Lyons},
  \bibinfo{title}{Neural SDEs as Infinite-Dimensional GANs}, in:
  \bibinfo{booktitle}{Proceedings of the 38th International Conference on
  Machine Learning}, vol. \bibinfo{volume}{139}, \bibinfo{pages}{5453--5463},
  \bibinfo{year}{2021}{\natexlab{b}}.

\bibitem[{Ye et~al.(2022)Ye, Li, Ye, and Zhang}]{YE2022118}
\bibinfo{author}{R.~Ye}, \bibinfo{author}{X.~Li}, \bibinfo{author}{Y.~Ye},
  \bibinfo{author}{B.~Zhang}, \bibinfo{title}{DynamicNet: A time-variant ODE
  network for multi-step wind speed prediction}, \bibinfo{journal}{Neural
  Networks} \bibinfo{volume}{152} (\bibinfo{year}{2022})
  \bibinfo{pages}{118--139}.

\bibitem[{Chu et~al.(2024{\natexlab{a}})Chu, Ma, Li, and Chen}]{CHU2024106549}
\bibinfo{author}{Z.~Chu}, \bibinfo{author}{W.~Ma}, \bibinfo{author}{M.~Li},
  \bibinfo{author}{H.~Chen}, \bibinfo{title}{Adaptive Decision Spatio-temporal
  neural ODE for traffic flow forecasting with Multi-Kernel Temporal Dynamic
  Dilation Convolution}, \bibinfo{journal}{Neural Networks}
  \bibinfo{volume}{179} (\bibinfo{year}{2024}{\natexlab{a}})
  \bibinfo{pages}{106549}.

\bibitem[{Chen et~al.(2018)Chen, Rubanova, Bettencourt, and
  Duvenaud}]{NEURIPS2018_69386f6b}
\bibinfo{author}{R.~T.~Q. Chen}, \bibinfo{author}{Y.~Rubanova},
  \bibinfo{author}{J.~Bettencourt}, \bibinfo{author}{D.~K. Duvenaud},
  \bibinfo{title}{Neural Ordinary Differential Equations}, in:
  \bibinfo{booktitle}{Advances in Neural Information Processing Systems},
  vol.~\bibinfo{volume}{31}, \bibinfo{year}{2018}.

\bibitem[{Rodriguez et~al.(2022)Rodriguez, Ames, and
  Yue}]{pmlr-v162-rodriguez22a}
\bibinfo{author}{I.~D.~J. Rodriguez}, \bibinfo{author}{A.~Ames},
  \bibinfo{author}{Y.~Yue}, \bibinfo{title}{{L}ya{N}et: A {L}yapunov Framework
  for Training Neural {ODE}s}, in: \bibinfo{booktitle}{Proceedings of the 39th
  International Conference on Machine Learning}, vol. \bibinfo{volume}{162},
  \bibinfo{pages}{18687--18703}, \bibinfo{year}{2022}.

\bibitem[{Kolter and Manek(2019)}]{NEURIPS2019_0a4bbced}
\bibinfo{author}{J.~Z. Kolter}, \bibinfo{author}{G.~Manek},
  \bibinfo{title}{Learning Stable Deep Dynamics Models}, in:
  \bibinfo{booktitle}{Advances in Neural Information Processing Systems},
  vol.~\bibinfo{volume}{32}, \bibinfo{year}{2019}.

\bibitem[{Kang et~al.(2021)Kang, Song, Ding, and Tay}]{NEURIPS2021_7d5430cf}
\bibinfo{author}{Q.~Kang}, \bibinfo{author}{Y.~Song},
  \bibinfo{author}{Q.~Ding}, \bibinfo{author}{W.~P. Tay},
  \bibinfo{title}{Stable Neural ODE with Lyapunov-Stable Equilibrium Points for
  Defending Against Adversarial Attacks}, in: \bibinfo{booktitle}{Advances in
  Neural Information Processing Systems}, vol.~\bibinfo{volume}{34},
  \bibinfo{pages}{14925--14937}, \bibinfo{year}{2021}.

\bibitem[{Zhao et~al.(2023)Zhao, Kang, Song, She, Wang, and
  Tay}]{NEURIPS2023_0a443a00}
\bibinfo{author}{K.~Zhao}, \bibinfo{author}{Q.~Kang},
  \bibinfo{author}{Y.~Song}, \bibinfo{author}{R.~She},
  \bibinfo{author}{S.~Wang}, \bibinfo{author}{W.~P. Tay},
  \bibinfo{title}{Adversarial Robustness in Graph Neural Networks: A
  Hamiltonian Approach}, in: \bibinfo{booktitle}{Advances in Neural Information
  Processing Systems}, vol.~\bibinfo{volume}{36}, \bibinfo{pages}{3338--3361},
  \bibinfo{year}{2023}.

\bibitem[{Cui et~al.(2023)Cui, Zhang, Chu, Hu, and Li}]{CUI2023576}
\bibinfo{author}{W.~Cui}, \bibinfo{author}{H.~Zhang}, \bibinfo{author}{H.~Chu},
  \bibinfo{author}{P.~Hu}, \bibinfo{author}{Y.~Li}, \bibinfo{title}{On
  robustness of neural ODEs image classifiers}, \bibinfo{journal}{Information
  Sciences} \bibinfo{volume}{632} (\bibinfo{year}{2023})
  \bibinfo{pages}{576--593}.

\bibitem[{Chu et~al.(2024{\natexlab{b}})Chu, Wei, Lu, and
  Zhao}]{https://doi.org/10.1049/cvi2.12248}
\bibinfo{author}{H.~Chu}, \bibinfo{author}{S.~Wei}, \bibinfo{author}{Q.~Lu},
  \bibinfo{author}{Y.~Zhao}, \bibinfo{title}{Improving neural ordinary
  differential equations via knowledge distillation}, \bibinfo{journal}{IET
  Computer Vision} \bibinfo{volume}{18}~(\bibinfo{number}{2})
  (\bibinfo{year}{2024}{\natexlab{b}}) \bibinfo{pages}{304--314}.

\bibitem[{Fazlyab et~al.(2019)Fazlyab, Robey, Hassani, Morari, and
  Pappas}]{NEURIPS2019_95e1533e}
\bibinfo{author}{M.~Fazlyab}, \bibinfo{author}{A.~Robey},
  \bibinfo{author}{H.~Hassani}, \bibinfo{author}{M.~Morari},
  \bibinfo{author}{G.~Pappas}, \bibinfo{title}{Efficient and Accurate
  Estimation of Lipschitz Constants for Deep Neural Networks}, in:
  \bibinfo{booktitle}{Advances in Neural Information Processing Systems},
  vol.~\bibinfo{volume}{32}, \bibinfo{year}{2019}.

\bibitem[{Latorre et~al.(2020)Latorre, Rolland, and
  Cevher}]{Latorre2020Lipschitz}
\bibinfo{author}{F.~Latorre}, \bibinfo{author}{P.~Rolland},
  \bibinfo{author}{V.~Cevher}, \bibinfo{title}{Lipschitz constant estimation of
  Neural Networks via sparse polynomial optimization}, in:
  \bibinfo{booktitle}{International Conference on Learning Representations},
  \bibinfo{year}{2020}.

\bibitem[{Khalil(2002)}]{Khalil:1173048}
\bibinfo{author}{H.~K. Khalil}, \bibinfo{title}{{Nonlinear systems; 3rd ed.}},
  \bibinfo{publisher}{Prentice-Hall}, \bibinfo{address}{Upper Saddle River,
  NJ}, \bibinfo{year}{2002}.

\bibitem[{Garg et~al.(2022)Garg, Arabi, and Panagou}]{GARG2022110314}
\bibinfo{author}{K.~Garg}, \bibinfo{author}{E.~Arabi},
  \bibinfo{author}{D.~Panagou}, \bibinfo{title}{Fixed-time control under
  spatiotemporal and input constraints: A Quadratic Programming based
  approach}, \bibinfo{journal}{Automatica} \bibinfo{volume}{141}
  (\bibinfo{year}{2022}) \bibinfo{pages}{110314}.

\bibitem[{YAN et~al.(2020)YAN, DU, TAN, and FENG}]{YAN2020On}
\bibinfo{author}{H.~YAN}, \bibinfo{author}{J.~DU}, \bibinfo{author}{V.~TAN},
  \bibinfo{author}{J.~FENG}, \bibinfo{title}{On Robustness of Neural Ordinary
  Differential Equations}, in: \bibinfo{booktitle}{International Conference on
  Learning Representations}, \bibinfo{year}{2020}.

\bibitem[{Zakwan et~al.(2023)Zakwan, Xu, and Ferrari-Trecate}]{9809979}
\bibinfo{author}{M.~Zakwan}, \bibinfo{author}{L.~Xu},
  \bibinfo{author}{G.~Ferrari-Trecate}, \bibinfo{title}{Robust Classification
  Using Contractive Hamiltonian Neural ODEs}, \bibinfo{journal}{IEEE Control
  Systems Letters} \bibinfo{volume}{7} (\bibinfo{year}{2023})
  \bibinfo{pages}{145--150}.

\bibitem[{Garg(2021)}]{garg2021advances}
\bibinfo{author}{K.~Garg}, \bibinfo{title}{Advances in the theory of fixed-time
  stability with applications in constrained control and optimization}, Ph.D.
  thesis, \bibinfo{school}{University of Michigan, Horace H. Rackham School of
  Graduate Studies}, \bibinfo{year}{2021}.

\bibitem[{Wang et~al.(2024)Wang, Ke, Hu, and Wu}]{WANG2024106087}
\bibinfo{author}{R.~Wang}, \bibinfo{author}{H.~Ke}, \bibinfo{author}{M.~Hu},
  \bibinfo{author}{W.~Wu}, \bibinfo{title}{Adversarially robust neural networks
  with feature uncertainty learning and label embedding},
  \bibinfo{journal}{Neural Networks} \bibinfo{volume}{172}
  (\bibinfo{year}{2024}) \bibinfo{pages}{106087}.

\bibitem[{Li et~al.(2023)Li, Zhang, Cao, and Tan}]{LI2023177}
\bibinfo{author}{J.~Li}, \bibinfo{author}{S.~Zhang}, \bibinfo{author}{J.~Cao},
  \bibinfo{author}{M.~Tan}, \bibinfo{title}{Learning defense transformations
  for counterattacking adversarial examples}, \bibinfo{journal}{Neural
  Networks} \bibinfo{volume}{164} (\bibinfo{year}{2023})
  \bibinfo{pages}{177--185}.

\bibitem[{Xiao(2017)}]{xiao2017fashion}
\bibinfo{author}{H.~Xiao}, \bibinfo{title}{Fashion-mnist: a novel image dataset
  for benchmarking machine learning algorithms}, \bibinfo{journal}{arXiv
  preprint arXiv:1708.07747} .

\bibitem[{Krizhevsky et~al.(2009)}]{krizhevsky2009learning}
\bibinfo{author}{A.~Krizhevsky}, et~al., \bibinfo{title}{Learning multiple
  layers of features from tiny images} .

\bibitem[{Liu et~al.(2021)Liu, Bernstein, Meister, and Yue}]{pmlr-v139-liu21c}
\bibinfo{author}{Y.~Liu}, \bibinfo{author}{J.~Bernstein},
  \bibinfo{author}{M.~Meister}, \bibinfo{author}{Y.~Yue},
  \bibinfo{title}{Learning by Turning: Neural Architecture Aware Optimisation},
  in: \bibinfo{booktitle}{Proceedings of the 38th International Conference on
  Machine Learning}, vol. \bibinfo{volume}{139}, \bibinfo{pages}{6748--6758},
  \bibinfo{year}{2021}.

\bibitem[{Hendrycks and Dietterich(2019)}]{hendrycks2018benchmarking}
\bibinfo{author}{D.~Hendrycks}, \bibinfo{author}{T.~Dietterich},
  \bibinfo{title}{Benchmarking Neural Network Robustness to Common Corruptions
  and Perturbations}, in: \bibinfo{booktitle}{International Conference on
  Learning Representations}, \bibinfo{year}{2019}.

\bibitem[{Goodfellow et~al.(2014)Goodfellow, Shlens, and
  Szegedy}]{goodfellow2014explaining}
\bibinfo{author}{I.~J. Goodfellow}, \bibinfo{author}{J.~Shlens},
  \bibinfo{author}{C.~Szegedy}, \bibinfo{title}{Explaining and harnessing
  adversarial examples}, \bibinfo{journal}{arXiv preprint arXiv:1412.6572} .

\bibitem[{Kurakin et~al.(2018)Kurakin, Goodfellow, and
  Bengio}]{kurakin2018adversarial}
\bibinfo{author}{A.~Kurakin}, \bibinfo{author}{I.~J. Goodfellow},
  \bibinfo{author}{S.~Bengio}, \bibinfo{title}{Adversarial examples in the
  physical world}, \bibinfo{publisher}{Chapman and Hall/CRC},
  \bibinfo{year}{2018}.

\bibitem[{Madry et~al.(2018)Madry, Makelov, Schmidt, Tsipras, and
  Vladu}]{madry2018towards}
\bibinfo{author}{A.~Madry}, \bibinfo{author}{A.~Makelov},
  \bibinfo{author}{L.~Schmidt}, \bibinfo{author}{D.~Tsipras},
  \bibinfo{author}{A.~Vladu}, \bibinfo{title}{Towards Deep Learning Models
  Resistant to Adversarial Attacks}, in: \bibinfo{booktitle}{International
  Conference on Learning Representations}, \bibinfo{year}{2018}.

\bibitem[{Croce and Hein(2020)}]{pmlr-v119-croce20b}
\bibinfo{author}{F.~Croce}, \bibinfo{author}{M.~Hein}, \bibinfo{title}{Reliable
  evaluation of adversarial robustness with an ensemble of diverse
  parameter-free attacks}, in: \bibinfo{booktitle}{Proceedings of the 37th
  International Conference on Machine Learning}, vol. \bibinfo{volume}{119},
  \bibinfo{pages}{2206--2216}, \bibinfo{year}{2020}.

\bibitem[{He et~al.(2016)He, Zhang, Ren, and Sun}]{He_2016_CVPR}
\bibinfo{author}{K.~He}, \bibinfo{author}{X.~Zhang}, \bibinfo{author}{S.~Ren},
  \bibinfo{author}{J.~Sun}, \bibinfo{title}{Deep Residual Learning for Image
  Recognition}, in: \bibinfo{booktitle}{Proceedings of the IEEE Conference on
  Computer Vision and Pattern Recognition (CVPR)}, \bibinfo{pages}{770--778},
  \bibinfo{year}{2016}.

\end{thebibliography}
\onecolumn
\appendix

\section{Theoretical Proofs}
\label{sec:theoretical_proofs}

\subsection{Proof of Theorem \ref{th:1}}
\label{sec:A1}
\begin{proof}
	We begin by recalling that
\begin{equation}\label{eq:A_1}
\mathcal{V}(x, y, h, t) := \max \left\{ 0, \frac{\partial V}{\partial h}f(h, x_c, t;\theta_f) + \alpha_1 V(x)^{\gamma_1} + \alpha_2 V(x)^{\gamma_2}\right\}.
\end{equation}
Clearly, $\mathcal{V} \geq 0$. We can see that $\mathcal{V}$ is continuous for any $x$, $y$, $h$, and $t$. In fact, the continuity of $\mathcal{V}$ can be attributed to its nature as the maximum of two continuous functions. Now we rearrange the integral terms in Equation \eqref{eq:in_loss} in the following form
\begin{align}
	\mathcal{L}_{lyp}(\theta) := { \mathbb{E}_{(x,y) \sim D}}\left[ {\int_0^1 {{\cal V}\left( {x,y,{h_\theta }\left( t  \right),t } \right){\mkern 1mu} dt} } \right]  \label{eq:A_3},
\end{align}
where ${h_\theta }(\cdot)$ is a solution of Equations \eqref{eq:2.1}-\eqref{eq:2.3}. 

We first claim that, when $\mathcal{L}_{lyp}(\theta) = 0$, $\mathcal{V}$ equals zero for almost everywhere  \( (x, y) \sim D \) and any $t \in [0,1]$.   Assume to the contrary that there exist $t^* \in [0,1]$ and non-zero measurement set  $\tilde D \in D$ such that $\mathcal{V} > 0$, i.e., there is a constant $\varepsilon > 0$ making $\mathcal{V} \ge \varepsilon$. Since $\mathcal{V}$ is continuous, there exists a $\delta$-neighborhood $B \in [0,1]$ such that $\mathcal{V} \ge \varepsilon$ for any $(x,y) \in \tilde{D}$ and $t \in B$. Thus, one has 
\begin{align}
	\mathcal{L}_{lyp}(\theta) =& { \mathbb{E}_{(x,y) \sim D}}\left[ {\int_0^1 {{\cal V}\left( {x,y,{h_\theta }\left( t  \right),t } \right){\mkern 1mu} dt} } \right] \\
	\ge & { \mathbb{E}_{(x,y) \sim \tilde{D}}}\left[ {\int_B {{\cal V}\left( {x,y,{h_\theta }\left( t  \right),t } \right){\mkern 1mu} dt} } \right] \\
	\ge &  \varepsilon { \mathbb{E}_{(x,y) \sim \tilde{D}}}\left[ {\int_B {dt} } \right]>0, 
\end{align}
which contradicts the assumption $\mathcal{L}(\theta)=0$, and so $\mathcal{V}=0$ for almost everywhere  \( (x, y) \sim D \) and any $t \in [0,1]$. 

Next we show the second result. By the definition of $\mathcal{V}$ and the first result, we have the following inequality 
\begin{equation}
\mathop V\limits^{\cdot} ({h_\theta }\left( t \right))\le -{\alpha _1}V{({h_\theta }\left( t \right))^{{\gamma _1}}} - {\alpha _2}V{({h_\theta }\left( t \right))^{{\gamma _2}}},
\end{equation}
where \(\alpha_1, \alpha_2 > 0\), \(\gamma_1 = 1 + \frac{1}{\mu}\), and \(\gamma_2 = 1 - \frac{1}{\mu}\) with \(\mu > 1\). 
It follows from Lemma \ref{lem:2.1} that $\sup \{ T \geq 0 : h_x(t) = h^* \text{ for all } t \geq T \} \le \frac{\mu \pi}{2\sqrt{\alpha_1 \alpha_2}}$, which is desired.
\end{proof}

\subsection{Proof of Proposition \ref{propos:1} } \label{sec:A2}

\begin{proof}
	Firstly, for $1 < \mu < 2$,  it holds that
	\begin{equation}
		I = \int_{{V_0}}^{\bar V} {\frac{{dV}}{{-{\alpha _1}{V^{{\gamma _1}}} - {\alpha _2}{V^{{\gamma _2}}} + {\delta}}}}  \le \int_{{V_0}}^{\bar V} {\frac{{dV}}{{-{\alpha _1}{V^{{\gamma _1}}} - {\alpha _2}{V^{{\gamma _2}}} + {\delta}{v^{1-\mu}}{V^{1 - \frac{1}{\mu }}}}}}
	\end{equation}
	due to the fact that ${{v^{1-\mu}}{V^{1 - \frac{1}{\mu }}}} > 1$ for all $V \in [\bar{V},V_{0}]$. Substituting  \( z = V^{\frac{1}{\mu}} \), we have $dz = \frac{1}{\mu }{V^{\frac{1}{\mu } - 1}}dV$, ${z_0} = V_0^{\frac{1}{\mu }}$, and ${\bar{z}} = \bar{V}^{\frac{1}{\mu }}$.  It follows that
	\begin{equation} \label{eq:A11}
	I \le  - \mu \int_{{z_0}}^{\bar z} {\frac{{{z^{\mu  - 1}}}}{{{\alpha _1}{z^{\mu  + 1}} + \left( {{\alpha _2} - \delta {v^{1 - \mu }}} \right){z^{\mu  - 1}}}}dz}  = - \frac{{\mu }}{{{\alpha _1}}}\int_{{z_0}}^{\bar z} {\frac{{dz}}{{\left( {z - v} \right)\left( {z + v} \right)}}}.
	\end{equation}
	Clearly, ${\left( {z - v} \right)\left( {z + v} \right)} >0 $ for all $z \in [\bar{z},z_{0}]$, which implies that the integral of the last equation in \eqref{eq:A11} exists. Thus, 
	\begin{equation}\label{eq:A12}
	 	I \le \frac{\mu }{{2{\alpha _1}v}}\left[ {\ln\left( {\left| {\frac{{\bar z + v}}{{\bar z - v}}} \right|} \right) + \ln\left( {\left| {\frac{{{z_0} - v}}{{{z_0} + v}}} \right|} \right)} \right].
	\end{equation}
	Taking ${z_0} = V_0^{\frac{1}{\mu }}$ and ${\bar{z}} = \bar{V}^{\frac{1}{\mu }}$ into above equation yields the desired result.
	
	Secondly, for $2 \le \mu < 3$,  it holds that
	\begin{equation}
		I = \int_{{V_0}}^{\bar V} {\frac{{dV}}{{-{\alpha _1}{V^{{\gamma _1}}} - {\alpha _2}{V^{{\gamma _2}}} + {\delta}}}}  \le \int_{{V_0}}^{\bar V} {\frac{{dV}}{{-{\alpha _1}{V^{{\gamma _1}}} - {\alpha _2}{V^{{\gamma _2}}} + {\delta}{v^{2 - \mu }}{V^{1 - \frac{2}{\mu }}}}}}
	\end{equation}
	due to the fact that ${{v^{2-\mu}}{V^{1 - \frac{2}{\mu }}}} > 1$ for all $V \in [\bar{V},V_{0}]$. Similarly, for \( z = V^{\frac{1}{\mu}} \), we have
	\begin{equation} \label{eq:A15}
		I \le  - \mu \int_{{z_0}}^{\tilde z} {\frac{{{z^{\mu  - 1}}}}{{{\alpha _1}{z^{\mu  + 1}} + {\alpha _2}{z^{\mu  - 1}} - \delta {v^{2 - \mu }}{z^{\mu  - 2}}}}} dz =  - \mu \int_{{z_0}}^{\tilde z} {\frac{z}{{\left( {z - v} \right)\left( {{\alpha _1}{z^2} + {\alpha _1}vz + {\alpha _1}{v^2} + {\alpha _2}} \right)}}} dz.
	\end{equation}
	Since $a_{1}, a_{2}, v >0$, one has ${\left( {z - v} \right)\left( {{\alpha _1}{z^2} + {\alpha _1}vz + {\alpha _1}{v^2} + {\alpha _2}} \right)}>0$ for all $z \in [\bar{z},z_{0}]$. Thus,  the integral of the last equation in \eqref{eq:A15} exists. Find \(A\), \(B\), and \(C\) such that
	\begin{equation}
	\frac{z}{{\left( {z - v} \right)\left( {{\alpha _1}{z^2} + {\alpha _1}vz + {\alpha _1}{v^2} + {\alpha _2}} \right)}} = \frac{A}{{z - v}} + \frac{{Bz + C}}{{{\alpha _1}{z^2} + {\alpha _1}vz + {\alpha _1}{v^2} + {\alpha _2}}}, 
	\end{equation}
	which is equivalent to 
	\begin{equation}
		z = (A{\alpha _1} + B){z^2} + (A{\alpha _1}v - Bv + C)z + (A{\alpha _1}{v^2} + A{\alpha _2} - Cv).
	\end{equation}
Thus, we have
	\begin{equation}
		\left\{ {\begin{array}{*{20}{l}}
				{A{\alpha _1} + B = 0},\\
				{A{\alpha _1}v - Bv + C = 1},\\
				{A{\alpha _1}{v^2} + A{\alpha _2} - Cv = 0}.
		\end{array}} \right. 
	\end{equation}
Since $\alpha_1,\alpha_2,v >0$, one has
	\begin{equation}
A = \frac{v}{{3{\alpha_1}{v^2} + {\alpha_2}}},\quad B = \frac{{ - {\alpha_1}v}}{{3{\alpha_1}{v^2} + {\alpha_2}}},\quad C = \frac{{{\alpha_1}{v^2} + {\alpha_2}}}{{3{\alpha_1}{v^2} + {\alpha_2}}}.
	\end{equation}
It follows from Equation \eqref{eq:A15} that 
	\begin{align}
	& - \mu \int_{{z_0}}^{\bar z} {\frac{z}{{\left( {z - v} \right)\left( {{\alpha _1}{z^2} + {\alpha _1}vz + {\alpha _1}{v^2} + {\alpha _2}} \right)}}dz} \nonumber \\
	 =&\frac{{ - \mu }}{{3{\alpha _1}{v^2} + {\alpha _2}}}\left( {\int_{{z_0}}^{\bar z} {\frac{v}{{z - v}}dz}  - \int_{{z_0}}^{\bar z} {\frac{{vz - {v^2} - \frac{{{\alpha _2}}}{{{\alpha _1}}}}}{{{z^2} + vz + {v^2} + \frac{{{\alpha _2}}}{{{\alpha _1}}}}}dz} } \right) \nonumber \\
	=&\frac{{ - \mu }}{{3{\alpha _1}{v^2} + {\alpha _2}}}\left[ {\left. {\left( {v\ln \left( {z - v} \right) - \frac{v}{2}\ln \left( {{z^2} + vz + {v^2} + \frac{{{\alpha _2}}}{{{\alpha _1}}}} \right)} \right)} \right|_{{z_0}}^{\bar z} + \frac{{3{v^2} + \frac{{2{\alpha _2}}}{{{\alpha _1}}}}}{{\sqrt {3{v^2} + \frac{{4{\alpha _2}}}{{{\alpha _1}}}} }}\left. {\arctan \left( {\frac{{2z + v}}{{\sqrt {3{v^2} + \frac{{4{\alpha _2}}}{{{\alpha _1}}}} }}} \right)} \right|_{{z_0}}^{\bar z}} \right] \nonumber \\
	 =&\frac{{\mu v}}{{6{\alpha _1}{v^2} + 2{\alpha _2}}}\left[ {\ln \left| {\frac{{{{\bar z}^2} + v\bar z + {v^2} + \frac{{{\alpha _2}}}{{{\alpha _1}}}}}{{{{\left( {\bar z - v} \right)}^2}}}} \right| + \ln \left| {\frac{{{{\left( {{z_0} - v} \right)}^2}}}{{{z_0}^2 + v{z_0} + {v^2} + \frac{{{\alpha _2}}}{{{\alpha _1}}}}}} \right| + \frac{{6v + \frac{{4{\alpha _2}}}{{{\alpha _1}v}}}}{{\sqrt {3{v^2} + \frac{{4{\alpha _2}}}{{{\alpha _1}}}} }}\left( {\arctan \left( {\frac{{2{z_0} + v}}{{\sqrt {3{v^2} + \frac{{4{\alpha _2}}}{{{\alpha _1}}}} }}} \right) - \arctan \left( {\frac{{2\bar z + v}}{{\sqrt {3{v^2} + \frac{{4{\alpha _2}}}{{{\alpha _1}}}} }}} \right)} \right)} \right]. \nonumber 
	\end{align}
	Taking ${z_0} = V_0^{\frac{1}{\mu }}$ and ${\bar{z}} = \bar{V}^{\frac{1}{\mu }}$ into above equation yields desired result.
	
	Finally, for $\mu \ge 3$, it holds that
	\begin{equation}
		I = \int_{{V_0}}^{\bar V} {\frac{{dV}}{{-{\alpha _1}{V^{{\gamma_1}}} - {\alpha_2}{V^{{\gamma_2}}} + {\delta _1}}}}  \le \int_{{V_0}}^{\bar V} {\frac{{dV}}{{-{\alpha_1}{V^{{\gamma_1}}} - {\alpha_2}{V^{{\gamma_2}}} + {\delta}{v^{3 - \mu }}{V^{1 - \frac{3}{\mu }}}}}}
	\end{equation}
due to the fact that ${{v^{3-\mu}}{V^{1 - \frac{3}{\mu }}}} > 1$ for all $V \in [\bar{V},V_{0}]$. Similarly, for \( z = V^{\frac{1}{\mu}} \), we have
	\begin{equation}\label{eq:A22}
		I \le  - \mu \int_{{z_0}}^{\bar z} {\frac{{{z^{\mu  - 1}}}}{{{\alpha _1}{z^{\mu  + 1}} + {\alpha _2}{z^{\mu  - 1}} - \delta {v^{3 - \mu }}{z^{\mu  - 3}}}}} dz = \frac{{ - \mu }}{{2{\alpha _1}{v^2} + {\alpha_2}}}\int_{{z_0}}^{\bar z} {\frac{{{v^2}}}{{{z^2} - {v^2}}} + \frac{{{v^2} + \frac{{{\alpha _2}}}{{{\alpha _1}}}}}{{{z^2} + {v^2} + \frac{{{\alpha _2}}}{{{\alpha _1}}}}}} dz.
	\end{equation}
Obviously, the integral of the last equation in \eqref{eq:A22} exists. It follows that
	\begin{equation}
	I \le \frac{\mu }{{2{\alpha _1}{v^2} + {\alpha _2}}}\left[ {\frac{v}{2}\left( {\ln \left| {\frac{{\bar z + v}}{{\bar z - v}}} \right| + \ln \left| {\frac{{{z_0} - v}}{{{z_0} + v}}} \right|} \right) + {k_5}\left( {\arctan \left( {\frac{{{z_0}}}{{{k_5}}}} \right) - \arctan \left( {\frac{{\bar z}}{{{k_5}}}} \right)} \right)} \right], \quad {k_5} = \sqrt {{v^2} + \frac{{{\alpha _2}}}{{{\alpha _1}}}}.
	\end{equation}
Taking ${z_0} = V_0^{\frac{1}{\mu }}$ and ${\bar{z}} = \bar{V}^{\frac{1}{\mu }}$ into above equation yields desired result.
\end{proof}

\subsection{Proof of Theorem \ref{the:2}}
\label{sec:A3}
\begin{proof}
	We will use the notation for the time derivative of $V_y$ as follows
\begin{equation}
	{\dot V_y}\left( {h\left( t \right),x,t} \right) = \frac{d}{{dt}}{V_y}\left( {h\left( t \right)} \right) = \left. {\frac{{\partial {V_y}}}{{\partial h}}} \right|_{h = h(t)}^{\top} f\left( {h\left( t \right),\phi \left( {x;{\theta _\phi }} \right),t;{\theta _f}} \right). \nonumber	\end{equation}
Thus it follows that	
\begin{align}
	&\dot{V}_y(h, x + \rho, t)  \nonumber \\
	=& \dot{V}_y(h, x, t) + \dot{V}_y(h, x + \rho, t) - \dot{V}_y(h, x, t)    \nonumber\\
	\leq& \dot{V}_y(h, x, t) + \left| \dot{V}_y(h, x + \rho, t) - \dot{V}_y(h, x, t) \right|   \nonumber\\
	=& {\dot V_y}(h,x,t) + \left| {{{\frac{{\partial {V_y}}}{{\partial h}}}^{\top}}f\left( {h,\phi \left( {x + \rho ;{\theta _\phi }} \right),t;{\theta _f}} \right) - {{\frac{{\partial {V_y}}}{{\partial h}}}^{\top}}f\left( {h,\phi \left( {x;{\theta _\phi }} \right),t;{\theta _f}} \right)} \right|  \nonumber\\
	\leq& \dot{V}_y(h, x, t) + \underbrace {{L_V}{L_f}{L_{\phi}}}_L\| \rho \| \quad   \nonumber\\
	\leq& \dot{V}_y(h, x, t) + L\rho \quad   \nonumber\\
	\leq& -\alpha_1 V_y(h)^{ 1+ \frac{1}{\mu}} - \alpha_2 V_y(h)^{1 - \frac{1}{\mu}} + L\rho.
\end{align}
In order to use the Comparison Principle of differential equations, consider an auxiliary differential equation given by 
\begin{equation}
	\dot{y} = -\alpha_1 y^{1 + \frac{1}{\mu}} - \alpha _2 y^{1 - \frac{1}{\mu}}+L\rho \label{eq:A_14}.
\end{equation}
Let  $v>0$ such that $a_{1}v^{\mu+1}+a_2v^{\mu-1}=L \rho$. Then it follows that $V_0 \ge \bar{V}=(\gamma v)^{\mu}$ with $\gamma >1$.  By using Proposition \ref{propos:1},  we can conclude that 
\begin{enumerate}
	\item  Since $v>0$, one has $V_0^{\frac{1}{\mu }} - v < V_0^{\frac{1}{\mu }} + v$. For $1 < \mu < 2$,  we can obtain an estimate of the fixed time as follows
	\begin{equation}
		T = \frac{\mu }{{2{\alpha _1}v}}\ln\left( {\left| {\frac{{(1 + \gamma )v}}{{(\gamma  - 1)v}}} \right|} \right).
	\end{equation}
	\item Since $v>0$, one has ${\left( {{V_0}^{\frac{1}{\mu }} - v} \right)^2} < {V_0}^{\frac{2}{\mu }} + v{V_0}^{\frac{1}{\mu }} + {v^2} + \frac{{{\alpha_2}}}{{{\alpha_1}}}$. By the fact that $\arctan \left( \cdot \right) < \frac{\pi }{2}$, for $2 \le \mu < 3$, we can obtain an estimate of the fixed time as follows
	\begin{equation}
		T = \frac{{\mu v}}{{6{\alpha _1}{v^2} + 2{\alpha _2}}}\left[ {\ln \left| {\frac{{\left( {{\gamma ^2} + \gamma  + 1} \right){v^2} + \frac{{{\alpha _2}}}{{{\alpha _1}}}}}{{{{\left( {\left( {\gamma  - 1} \right)v} \right)}^2}}}} \right| + \frac{{6v + \frac{{4{\alpha _2}}}{{{\alpha _1}v}}}}{{\sqrt {3{v^2} + \frac{{4{\alpha _2}}}{{{\alpha _1}}}} }}\left( {\frac{\pi }{2} - \arctan \frac{{\left( {2\gamma  + 1} \right)v}}{{\sqrt {3{v^2} + \frac{{4{\alpha _2}}}{{{\alpha _1}}}} }}} \right)} \right] .
	\end{equation}
	
	\item  Since $v>0$, one has $V_0^{\frac{1}{\mu }} - v < V_0^{\frac{1}{\mu }} + v$. By the fact that $\arctan \left( \cdot \right) < \frac{\pi }{2}$,  for $\mu \ge 3$, we can obtain an estimate of the fixed time as follows
	\begin{equation}
		T = \frac{\mu }{{2{\alpha _1}{v^2} + {\alpha _2}}}\left[ {\frac{v}{2}\left( {\ln \left| {\frac{{1 + \gamma }}{{1 - \gamma }}} \right|} \right) + {k_5}\left( {\frac{\pi }{2} - \arctan \left( {\frac{{\gamma v}}{{{k_5}}}} \right)} \right)} \right],
	\end{equation}
	where $k_5=\sqrt {{v^2} + \frac{{{\alpha _2}}}{{{\alpha _1}}}}$.
\end{enumerate}
In summary, when $t \ge T$, for any $V_y(h_{\tilde{x}}(0)) \ge (\gamma v)^{\mu}$, one has $V_y(h_{\tilde{x}}(t))=\left\| {h{}_{\tilde{x}}\left( t \right) - {h^*}} \right\|_2^2 \le (\gamma v)^{\mu}$ with $\gamma >1$. Thus, we obtain $\left| {\tilde y(t) - y} \right| \le {\left( {\gamma v} \right)^\mu }{L_\psi }$ with $\tilde{y}(t)=\psi \left( {h_{\tilde{x}}\left( t \right);{\theta _\psi }} \right)$, which is desired.
\end{proof}

\newpage
\section{Experiments}
\label{sec:Experiments}
\begin{table}[th]
		\centering
	\caption{ Classification error and robustness to perturbations for the ablation study of the parameter $\eta_{2}$, where the adversarial attack perturbations $\epsilon=8/255$. ``Original" denotes the classification error on original data. The best result is in bold. ``Avg" denotes the average error. Unit: \%.} \label{table:ab3}
	\begin{tabular}{cccccccccc}
		\toprule
		Parameters                         & Orignal        & Gaussian       & Impulse        & Speckle        & Motion         & FGSM           & BIM            & APGD         & Avg            \\ \hline
		$\eta_{2}=0.4$                     & 17.31          & 45.81          & 55.17          & 22.45          & 25.66          & 37.94          & 19.73          & 40.06          & 33.02          \\
		$\eta_{2}=0.8$                     & 16.99          & 44.79          & 53.50          & 22.43          & 25.17          & 37.82          & 19.73          & 40.17          & 32.58          \\
		$\eta_{2}=1.2$                     & 16.91          & 44.27          & 53.32          & 22.56          & 24.91          & 37.80          & 19.81          & 39.89          & 32.43          \\
		$\eta_{2}=1.6$                     & 16.80          & 44.48          & 53.08          & 21.79          & \textbf{24.37} & 37.40          & 19.75          & 39.52          & 32.15          \\
		$\eta_{2}=2.0$ & \textbf{16.21} & 43.96          & \textbf{52.64} & \textbf{21.63} & 24.47          & \textbf{37.37} & \textbf{19.15} & \textbf{39.48} & \textbf{31.86} \\
		$\eta_{2}=2.4$ & 16.64          & \textbf{43.44} & 52.83          & 22.15          & 25.46          & 37.87          & 19.64          & 39.74          & 32.22          \\
		$\eta_{2}=2.8$ & 17.37          & 45.69          & 54.38          & 22.53          & 25.32          & 37.66          & 19.96          & 39.98          & 32.86 \\  \bottomrule             
	\end{tabular}
\end{table}
Following the settings of Section \ref{sec:ASP}, we investigated the impact of the inner learning rate $\eta_2$ on the classification performance and robustness of FxTS-Net, where the parameter $\eta_2$ varies between 0.4 and 2.8. Table \ref{table:ab3} shows a general decreasing trend in classification error as $\eta_2$ increases, with the smallest error of 16.21\% at $\eta_2 = 2.0$ on the original data. $\eta_2 = 2.0$ consistently showed strong performance for adversarial attacks and noise, producing the lowest average classification error (31.86\%). While larger values of $\eta_2$ slightly reduce robustness, increasing $\eta_2$ too much (e.g., $\eta_2 = 2.8$) tends to increase classification error, especially with Gaussian and impulse noise. The optimal value of $\eta_2$ is around $2.0$, which balances accuracy and robustness to various perturbations. 
\begin{table}[th]
		\centering
	\caption{ Classification error and robustness to perturbations for the ablation study of the parameter $\Gamma,n_{\delta }$, where the adversarial attack perturbations $\epsilon=8/255$. ``Original" denotes the classification error on original data. The best result is in bold. ``Avg" denotes the average error. Unit: \%.} \label{table:ab4}
	\begin{tabular}{lccccccccc}
		\toprule
		\multicolumn{1}{c}{Parameters} & Orignal        & Gaussian       & Impulse        & Speckle        & Motion         & FGSM           & BIM            & APGD       & AVG            \\ \hline
		$\Gamma=3,n_{\delta }=170$     & 16.79          & \textbf{43.79} & 53.00          & 22.13          & 25.23          & 37.45          & 19.71          & 39.59          & 32.21          \\
		$\Gamma=5,n_{\delta }=102$      & \textbf{16.21} & 43.96          & \textbf{52.64} & 21.63          & \textbf{24.47} & \textbf{37.37} & \textbf{19.15} & \textbf{39.48} & \textbf{31.86} \\
		$\Gamma=7,n_{\delta }=73$      & 16.69          & 44.70          & 53.27          & \textbf{21.41} & 24.64          & 38.32          & 19.94          & 40.31          & 32.41          \\
		$\Gamma=9,n_{\delta }=56$      & 16.73          & 44.60          & 53.28          & 22.46          & 24.69          & 38.13          & 19.71          & 40.26          & 32.48          \\
		$\Gamma=11,n_{\delta }=46$     & 17.00          & 45.26          & 53.70          & 22.04          & 24.83          & 38.12          & 19.95          & 40.45          & 32.67          \\
		$\Gamma=13,n_{\delta }=39$     & 17.03          & 45.54          & 54.20          & 22.29          & 25.80          & 38.32          & 20.04          & 40.39          & 32.95    \\  \bottomrule 
	\end{tabular}
\end{table}

We consider the effect of the trajectory sampling resolution $\Gamma$ and the number of perturbation samples $n_\delta$ on the model performance, where $\Gamma, n_\delta$ satisfies $\Gamma* n_\delta \le 512$. As shown in Table \ref{table:ab4}, the configuration $\Gamma=5$ and $n_\delta=102$ achieves the lowest classification error on both clean data and robustness. As $\Gamma$ increases beyond 5, the model's robustness to perturbations starts to degrade, as shown by the increasing average error, particularly for $\Gamma=13$. Similarly, decreasing the number of perturbation samples $n_\delta$ leads to poorer performance, likely due to the reduced diversity in the perturbation space that the model is exposed to during training. The combination of $\Gamma=5$ and $n_\delta=102$ strikes an optimal balance between model robustness and classification accuracy. Higher values of $\Gamma$ or lower values of $n_\delta$ lead to a deterioration in performance, indicating that careful tuning of both parameters is critical for robustness





\end{document}